\documentclass{amsart}
\usepackage{amscd}
\usepackage{amssymb}
\usepackage{graphics}
\usepackage{pst-tree}

\numberwithin{equation}{section}

\newtheorem{Thm}{Theorem}[section]
\newtheorem{Prop}[Thm]{Proposition}
\newtheorem{Lem}[Thm]{Lemma}
\newtheorem{Cor}[Thm]{Corollary}

\theoremstyle{remark}
\newtheorem{Rem}[Thm]{Remark}
\newtheorem*{Ack}{Acknowledgment}

\theoremstyle{definition}
\newtheorem{Def}[Thm]{Definition}
\newtheorem{Exa}[Thm]{Example}

\newtheorem{Conv}[Thm]{Convention}
\newtheorem{Conj}[Thm]{Conjecture}

\def\bconv{\begin{Conv}}
\def\econv{\end{Conv}}
\def\bprop{\begin{Prop}}
\def\eprop{\end{Prop}}
\def\bconj{\begin{Conj}}
\def\econj{\end{Conj}}
\def\bth{\begin{Thm}}
\def\eth{\end{Thm}}
\def\bcor{\begin{Cor}}
\def\ecor{\end{Cor}}
\def\blem{\begin{Lem}}
\def\elem{\end{Lem}}
\def\bdefiniz{\begin{Def}}
\def\edefiniz{\end{Def}}

\newcommand{\fullref}[1]{\ref{#1} on page~\pageref{#1}}

\DeclareMathOperator{\Hom}{Hom}

\newcommand{\lora}{{\longrightarrow}}

\newcommand\de{\partial}

\newcommand{\Z}{\mathbb{Z}}
\newcommand{\N}{\mathbb{N}}

\newcommand{\End}{\mathcal E\!nd}

\newcommand{\Det}{\mathrm{Det}}

\newcommand{\f}{
 \begin{pspicture}(0,0.1)(0.2,0.4)
 \psline[linewidth=1pt](0.1,0)(0.1,0.4)
 \end{pspicture}}
\newcommand{\da}{
 \begin{pspicture}(0,0.1)(0.2,0.4)
 \psline[linewidth=1pt, linestyle=dashed,dash=4pt 3pt](0.1,0)(0.1,0.4)
 \end{pspicture}}
\newcommand{\ff}{
 \begin{pspicture}(0,0)(0.2,0.2)
 \psline[linewidth=1pt](0.1,0)(0.1,0.2)
 \end{pspicture}}
\newcommand{\dda}{
 \begin{pspicture}(0,0)(0.2,0.2)
 \psline[linewidth=1pt, linestyle=dashed,dash=3pt 2pt](0.1,0)(0.1,0.2)
 \end{pspicture}}

\def\cmp#1,{{ Commun.\ Math.\ Phys.\ \bf #1},}
\def\jmp#1,{{ J.\ Math.\ Phys.\ \bf #1},}
\def\pl#1,{{ Phys.\ Lett.\ \bf #1},}
\def\npb#1,{{ Nucl.\ Phys.\ {\bf B #1}},}
\def\mpl#1,{{ Mod.\ Phys.\ Lett.\ \bf #1},}
\def\pr#1,{{ Phys.\ Rev.\ \bf #1},}
\def\prl#1,{{ Phys.\ Rev.\ Lett.\ \bf #1},}
\def\lmp#1,{{ Lett.\ Math.\ Phys.\ \bf #1},}
\def\jktr#1,{{ Jour.\ of Knot Theory and its Ramification \bf #1},}
\def\bams#1,{{ Bull.\ Amer.\ Math.\ Soc.\ \bf #1},}
\def\ja#1,{{ J.\ Algebra \bf #1},}
\def\ijm#1,{{ Illinois Journal of Mathematics \bf #1},}
\def\dmj#1,{{ Duke Math.\ J.\ \bf #1},}
\def\lnm#1,{{ Lecture Notes in Mathematics \bf #1},}
\def\mrl#1,{{ Math.\ Res.\ Lett.\ \bf #1},}

\begin{document}

\title{Homotopy Inner Products for Cyclic Operads}

\author[R.~Longoni]{Riccardo~Longoni}
\address{Dipartimento di Matematica ``G. Castelnuovo'',
Universit\`a di Roma ``La Sapienza'', Piazzale Aldo Moro, 2
I-00185 Roma, Italy}
\email{longoni@mat.uniroma1.it}

\author[T.~Tradler]{Thomas~Tradler}
\address{College of Technology of the City University
of New York, Department of Mathematics, 300 Jay Street, Brooklyn,
NY 11201, USA}
\email{ttradler@citytech.cuny.edu}

\begin{abstract}
We introduce the notion of homotopy inner products for any cyclic
quadratic Koszul operad $\mathcal O$, generalizing the construction
already known for the associative operad. This is done by defining a
colored operad $\widehat{\mathcal O}$, which describes modules over
$\mathcal O$ with invariant inner products. We show that
$\widehat{\mathcal O}$ satisfies Koszulness and identify algebras
over a resolution of $\widehat{\mathcal O}$ in terms of derivations
and module maps. An application to Poincar\'e duality on the chain
level of a suitable topological space is given.
\end{abstract}


\maketitle

\tableofcontents

In \cite{GeK}, the notion of cyclic operads and invariant inner
product for such operads was defined. A homotopy version of these
inner products for the associative operad was given in \cite{Tr}.
Thus, it is natural to ask for a generalization of this construction
applicable to any cyclic operad. Following Markl's paper \cite{M} on
homotopy algebras, it is clear that in order to understand the
homotopy version of some object, one needs to resolve the
corresponding operad. Thus, a criterion for a suitable definition for
these homotopy inner products is, that it should describe the data of
an algebra over a resolution of some operad describing those inner
products. It turns out, that the correct setting for performing this
resolution is that of colored operads.

Here is an outline of the paper. In Section~\ref{sec:mod} we will review
the notion of colored operad  and adapt to this case the analysis made
in~\cite{GK} for ordinary operads. In particular the resolution
$\mathcal O_\infty$ of a Koszul quadratic colored operad $\mathcal O$
is described, by means of the ``cobar dual''. (In this note we will
assume familiarity with the notions of operads and
cyclic operads; for a good introduction to these topics, we refer the
reader to \cite{Ad}, \cite{GK}, \cite{GeK} and \cite{MSS}.) Then we use
this language to define, for any quadratic Koszul operad $\mathcal O$, a
colored operad $\overline{\mathcal O}$ whose algebras describe
modules over $\mathcal O$. It is shown, that $\overline{\mathcal O}$
satisfies the Koszulness property (Proposition~\ref{O_bar_Koszul})
and hence the homotopy version $\overline{\mathcal O}_\infty$ of
$\overline{\mathcal O}$ can be explicitly constructed using the cobar
dual. Algebras over $\overline{\mathcal O}_\infty$ are called
homotopy $\mathcal O$-modules and can be easily described in terms of
derivations (Theorem~\ref{O_bar_algebras}).

In Section~\ref{sec:cy}, after recalling the notion of cyclic operads
$\mathcal O$, we introduce a colored operad $\widehat{\mathcal O}$
whose algebras are modules over $\mathcal O$ endowed with an
invariant inner product. The main result of the paper is that
Koszulness of $\mathcal O$ implies Koszulness of $\widehat{\mathcal O}$
(Theorem~\ref{O_hat_Koszul}) and hence $\widehat{\mathcal O}_\infty$
can be explicitly constructed as well. It turns out that algebras
over $\widehat{\mathcal O}_\infty$ have some $n$-to-$0$ operations,
called {\em homotopy inner products}, which can be described in terms
of derivations (Theorem~\ref{O_hat_algebras}). When the associative
operad is considered we obtain the $\infty$-inner products introduced
in \cite{Tr}. Intuitively, Koszulness of $\widehat{\mathcal O}$ means
that one obtains a new resolution of $\mathcal O$, where elements of
$\mathcal O(n)$ are interpreted as $(n+1)$-to-$0$ diagrams,
see Remark~\ref{Koszulness_idea}.
Finally, as an application, we will recall and slightly generalize the
results of \cite{TZ}. More specifically, we will show in
Proposition~\ref{prop:comm-pd} that Poincar\`e duality on the chain level
of a topological space gives rise to a homotopy inner product structure.

\begin{Ack}
We are grateful to Dennis Sullivan for many valuable suggestions and
illuminating discussions. We also thank Domenico Fiorenza, Martin Markl
and Scott Wilson for useful comments and remarks regarding this topic.
\end{Ack}

\section{Algebras and Modules over Operads}
\label{sec:mod}

In the first three subsections, we repeat some standard definitions
and fix notation about colored quadratic operads, and Koszulness. In
the last subsection, we identify homotopy algebras and homotopy
modules in terms of derivations.

\subsection{Colored Operads}

Let $k$ be a field of characteristic $0$. In this paper we will be
concerned only with operads and PROPs over the category
$\underline{\mathrm{Vect}}$ and $\underline{\mathrm{dgVect}}$ of
vector spaces and differential graded vector spaces over $k$,
although all definitions can be made for any symmetric tensor
category. We use the convention that differential graded spaces have
differentials of degree $+1$, and the arrow of the differential goes
from left to right.

Roughly speaking a colored PROP is a similar to a PROP, except that
instead of considering collections of vector spaces $\mathcal P(n,l)$
for diagrams with $n$ inputs and $l$ outputs, $n,l\in\N$, one
replaces $n$ and $l$ with sequences of finite length of ``colored
elements''. When the possible colors are $m$, we have an $m$-colored
PROPs. Clearly when $m=1$ we get the usual definition of PROP.

\begin{Def}[$m$-colored Operad]
Let $\underline{m}$ denote the discrete tensor category having as
objects (nonempty) finite sequences of elements of $m$ types, and as
tensor product the joint of two sequences. An {\em$m$-colored PROP}
is a functor $\mathcal P\colon (\underline{m})^{op}\times
(\underline{m}) \to \underline{\mathrm{Vect}}$ (or $\mathcal P\colon
(\underline{m})^{op}\times (\underline{m}) \to \underline{
\mathrm{dgVect}}$) together with natural transformations
$\circ_{A,B,C}\colon \mathcal P(A,B)\otimes \mathcal P(B,C) \to
\mathcal P(A,C)$ called composition maps, identity maps $j_A$ in
$\mathcal P(A,A)$, natural transformations $\otimes_{A,B,C,D}\colon
P(A,B)\otimes \mathcal P(C,D) \to \mathcal P(A\otimes C,B\otimes D)$
and distinguished elements $\sigma_{A,B} \in \mathcal P(A\otimes B,
B\otimes A)$ such that all the natural diagrams commute (see e.g.
\cite{Ad}).\\%
An {\em$m$-colored dioperad} is the ``tree-part'' of an $m$-colored
PROP (see also \cite{Gan}).\\%
An {\em$m$-colored operad} is given when an $m$-colored dioperad is
generated only by $n$-to-1 operations.\\%
We will use the convention that when talking about operads, it is
meant that the operad is over \underline{Vect}, unless otherwise
specified. Some of the associated object to an operad, like for
example the cobar dual, will then be operads over \underline{dgVect}.
\end{Def}

In this paper we will be mainly interested in 2-colored operads, so
that it is helpful to write a more explicit version of the definition
of 2-colored operad given above. By convention we will denote the
colors by the symbols $\f$ (``full'') and $\da$ (``dashed'').
Thus, a 2-colored operad attaches multiplication data to trees of the
form
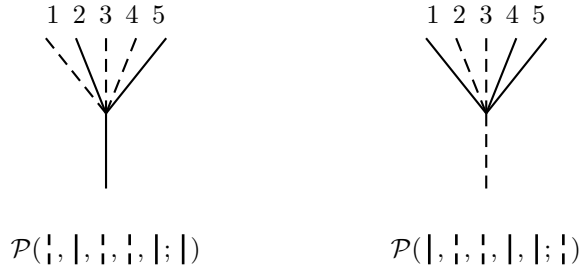
\[
\begin{pspicture}(0,0)(4,4)
 \psline[linestyle=dashed, arrowsize=0.1, arrowinset=0](2,2)(1.2,3)
 \psline[arrowsize=0.1, arrowinset=0](2,2)(1.6,3)
 \psline[linestyle=dashed, arrowsize=0.1, arrowinset=0](2,2)(2,3)
 \psline[linestyle=dashed, arrowsize=0.1, arrowinset=0](2,2)(2.4,3)
 \psline[arrowsize=0.1, arrowinset=0](2,2)(2.8,3)
 \psline[arrowsize=0.1, arrowinset=0](2,2)(2,1)
 \rput[b](2,0){$\mathcal P (\da,\f,\da,\da,\f;\f)$}
 \rput[b](2,3.2){$1\,\,\, 2\,\,\, 3\,\,\, 4\,\,\, 5$}
\end{pspicture}
\quad \quad \quad
\begin{pspicture}(0,0)(4,4)
 \psline[arrowsize=0.1, arrowinset=0](2,2)(1.2,3)
 \psline[linestyle=dashed, arrowsize=0.1, arrowinset=0](2,2)(1.6,3)
 \psline[linestyle=dashed, arrowsize=0.1, arrowinset=0](2,2)(2,3)
 \psline[arrowsize=0.1, arrowinset=0](2,2)(2.4,3)
 \psline[arrowsize=0.1, arrowinset=0](2,2)(2.8,3)
 \psline[linestyle=dashed, arrowsize=0.1, arrowinset=0](2,2)(2,1)
 \rput[b](2,0){$\mathcal P (\f,\da,\da,\f,\f;\da)$}
 \rput[b](2,3.2){$1\,\,\, 2\,\,\, 3\,\,\, 4\,\,\, 5$}
\end{pspicture}
\]
Unless otherwise stated, we use the convention to put inputs at the
top and evaluate towards the bottom output.\\ Explicitly, a 2-colored
operad $\mathcal P$ consists of a collection of $k$-vector spaces
$\mathcal P(\vec X;\f)$ and $\mathcal P(\vec X;\da)$ for every nonempty
finite sequence $\vec X$ of the symbols $\f$ and $\da$. In the tree
picture, the sequence $\vec X$ indicate which top (inputs) lines have
the ``full'' and which have the ``dashed'' color. The last argument of
$\mathcal P$ indicates the color of the bottom (output) line; see for
example the picture above.\\ We let the symbol $|\vec X|$ denote the
length of the sequence $\vec X$. The spaces $\bigoplus_{|\vec X|=n}
\mathcal P(\vec X;x)$, for $x\in \{\f,\da\}$, are required to be endowed
with an action of the symmetric group $S_n$, so that $\sigma\in S_n$
maps
\begin{equation}\label{equimap}
\sigma:\mathcal P(\vec X;x)\to \mathcal P(\sigma(\vec X);x)
\end{equation}
where the action of $\sigma$ on a sequence $\vec X$ of length $n$
simply permutes the order of the symbols $\f$ and $\da$ in the sequence.
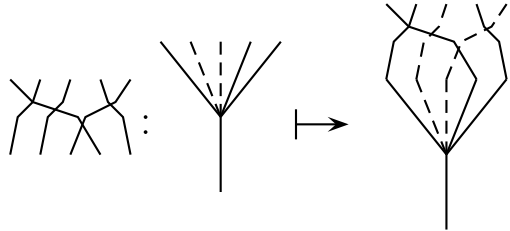
\[
\begin{pspicture}(1,1.5)(3,4.5)
 \psline(1.2,3)(1.3,3.5)(1.4,3.6)(1.5,3.7)(1.6,4)
 \psline(1.6,3)(1.7,3.5)(1.8,3.6)(1.9,3.7)(2,4)
 \psline(2,3)(2.2,3.5)(2.4,3.6)(2.6,3.7)(2.8,4)
 \psline(2.4,3)(2.1,3.5)(1.8,3.6)(1.5,3.7)(1.2,4)
 \psline(2.8,3)(2.7,3.5)(2.6,3.6)(2.5,3.7)(2.4,4)
 \psdot[dotsize=0.05](3,3.5) \psdot[dotsize=0.05](3,3.3)
\end{pspicture}
\begin{pspicture}(1,0)(3,4)
 \psline(2,2)(1.2,3)
 \psline[linestyle=dashed](2,2)(1.6,3)
 \psline[linestyle=dashed](2,2)(2,3)
 \psline(2,2)(2.4,3)
 \psline(2,2)(2.8,3)
 \psline(2,2)(2,1)
\end{pspicture}
\begin{pspicture}(0,0.5)(4,3.5)
 \psline(0,2.2)(0,2.6) \psline[arrowsize=0.2]{->}(0,2.4)(0.7,2.4)
 \psline(2,2)(1.2,3)  \psline(1.2,3)(1.3,3.5)(1.4,3.6)(1.5,3.7)(1.6,4)
 \psline[linestyle=dashed](2,2)(1.6,3)
 \psline[linestyle=dashed](1.6,3)(1.7,3.5)(1.8,3.6)(1.9,3.7)(2,4)
 \psline[linestyle=dashed](2,2)(2,3)
 \psline[linestyle=dashed](2,3)(2.2,3.5)(2.4,3.6)(2.6,3.7)(2.8,4)
 \psline(2,2)(2.4,3)  \psline(2.4,3)(2.1,3.5)(1.8,3.6)(1.5,3.7)(1.2,4)
 \psline(2,2)(2.8,3)  \psline(2.8,3)(2.7,3.5)(2.6,3.6)(2.5,3.7)(2.4,4)
 \psline(2,2)(2,1)
\end{pspicture}
\]
The definition of 2-colored operad includes composition maps
\begin{equation}
\label{diopmaps} \circ_{i}\colon \mathcal P(\vec X;x) \otimes
\mathcal P(\vec Y;y) \to \mathcal P(\vec Z;x),
\end{equation}
where the color $y$ has to be the same as the $i$th elements of
$\vec X$.
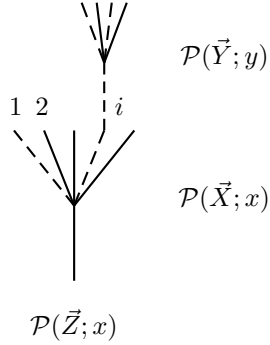
\[
\begin{pspicture}(0,0)(4,5)
 \psline[linestyle=dashed](2,2)(1.2,3)
 \psline(2,2)(1.6,3)
 \psline(2,2)(2,3)
 \psline[linestyle=dashed](2,2)(2.4,3)
 \psline(2,2)(2.8,3)
 \psline(2,2)(2,1)
 \rput[b](1.9,3.2){$1\,\,\, 2\quad\,\,\,\quad i$}
 \psline[linestyle=dashed](2.4,3)(2.4,3.9)
 \psline[linestyle=dashed](2.4,3.9)(2.1,4.7)
 \psline(2.4,3.9)(2.3,4.7)
 \psline[linestyle=dashed](2.4,3.9)(2.5,4.7)
 \psline(2.4,3.9)(2.7,4.7)
 \rput[b](4,1.9){$\mathcal P(\vec X;x)$}
 \rput[b](4,3.8){$\mathcal P(\vec Y;y)$}
 \rput[b](2,0.2){$\mathcal P(\vec Z;x)$}
\end{pspicture}
\]
These maps have to satisfy the usual associativity and equivariance
axioms stated for example in \cite{GK}.\\ \label{O(1)=k} In this
paper we want to make the assumption that, in the 2-colored case, the
only nontrivial spaces with one input are $\mathcal P(\f;\f)=k$ and
$\mathcal P(\da;\da)=k$, whereas in the 1-colored case, we only
assume $\mathcal P(\f;\f)=k$. For those spaces, we fix generators,
called {\em units}, which are required to satisfy the unit axioms as
stated in \cite{GK}.\\ An example of a 2-colored operad is given for
any two spaces $A$ and $M$ by the {\em endomorphism operad}
$\End^{A,M}$. Specifically, ${\End^{A,M}}(\vec X;x)$ is defined to be
the space of linear maps from tensor products of $A$ and $M$ to
either $A$ or $M$ according to $\vec X$ and $x\in \{\f,\da\}$, e.g.
${\End^{A,M}}(\da,\da,\f,\da,\f,\f;\da)= Hom(M\otimes M\otimes
A\otimes M\otimes A\otimes A, M)$. Here, $A$ corresponds to the
``full'' color and $M$ to the ``dashed'' color. The $S_n$ action is
defined by permuting the tensor products and the operad maps are
defined to be the composition of linear maps.

\begin{Def}[Algebras]\label{def_col_alg}
For any 2-colored operad $\mathcal P$, we define an algebra over
$\mathcal P$ to be a pair of $k$-vector spaces $(A,M)$ together with
a map of 2-colored operads $\mathcal P\lora \End^{A,M}$. By
definition this means that one has maps $\mathcal P(\vec X;x)\to
\End^{A,M}(\vec X;x)$, for every sequence $\vec X$ and every
$x\in\{\f,\da\}$ respecting equivariance, composition and units.
\end{Def}

An {\em operad} is a 1-colored operad. We will often simplify
notations by writing $\mathcal O (n)=\mathcal O (\f,\ldots,\f;\f)$.
Just as before, for a given space $A$, one has the {\em endomorphism
operad} $\End ^A$, defined by $\End ^A(n):=Hom(A^{\otimes n},A)$. We
say that the vector space $A$ is an algebra over an operad $\mathcal O$
if one has an operad map $\mathcal O \to \End^A$.

\begin{Rem}
All the above constructions can be carried over to the category of
differential graded vector spaces. In that case there exists a
differential acting on each $P(\vec X;x)$. In particular, given two
differential graded vector spaces $(A,d_A)$ and $(M,d_M)$, we have a
natural differential $\widetilde{\partial}^{A,M}$ on $\End^{A,M}(\vec
X;x)$ given by
\[
(\widetilde{\partial}^{A,M}(h))(x_1,\ldots,x_n) =
d_X(h(x_1,\ldots,x_n)) + \sum_{i=1}^n (-1)^{\epsilon_i}
h(x_1,\ldots,d_X(x_i),\ldots,x_n)
\]
where $d_X$ is either $d_A$ or $d_M$ depending whether it acts on an
element of $A$ or $M$ respectively, and $\epsilon_i = |h|+ |x_1|+
\cdots + |x_{i-1}|$\label{sgn_epsilon_i}.\\ Hence, an algebra over
the 2-colored operad $\mathcal P$ in the category of differential
graded spaces is a pair of differential graded spaces
$((A,d_A),(M,d_M))$ and a map $\mathcal P \to \End^{A,M}$ of
2-colored operad which respects the differential graded structure.\\
Clearly every vector space can be thought as a differential graded
vector space concentrated in degree zero with zero differential. In
this case it is $\widetilde{\partial}^{A,M}=0$, and we recover the
previous definitions.
\end{Rem}

\subsection{Cobar Dual for Colored Operads}

To every 2-colored operad, one can associate its 2-colored cobar dual
operad. This is defined using a sum of suitable trees. All the
definitions below, which are essentially an adaption of the ones in
\cite{GK}, are given for the 2-colored operads, although they can be
easily applied to the general case of $m$-colored (di)operads.

\begin{Def}[Tree]
A {\em tree} is a one-dimensional contractible complex consisting of
edges and vertices. The edges have one of two colors ``full'' or
``dashed''. The {\em valence} $|v|$ of the vertex $v$ is the number
of edges merging at $v$. Here we only consider trees with vertices of
valence $1$ or valence $\geq 3$. Edges that have a vertex of valence
$1$ in its boundary are called {\em external edges}. All other edges
are called {\em internal edges}. We require, that every vertex of
valence $\geq 3$ has exactly one ``output'' edge and $|v|-1$
``input'' edges. We denote by the symbol $\mathrm{In}(v)$ the set of
incoming edges to the vertex $v$. We can draw this in the plane by
putting the inputs on top and the output at the bottom of the vertex.
There is a unique edge which is external and is the output edge for
some vertex; this edge is called the {\em root}\/ or {\em output
edge}\/ of the tree. We call a vertex $v$ with $|v|\neq 1$ {\em of
type} $(\vec X;x)$ if the output edge of $v$ has color
$x\in\{\f,\da\}$ and there is a labeling of the incoming edges by
$\vec X$, so that an incoming edge is labeled by $\f$ if and only if
this edge has the ``full'' color. We call a tree $T$ {\em of type}
$(\vec X;x)$, if the output edge of the tree has color
$x\in\{\f,\da\}$ and there is a labeling of the incoming edges by
$\vec X$, so that an edge is labeled by $\f$ if and only if this edge
has the ``full'' color. Note that $\vec X=(x_1,...,x_n)$, $x_i\in
\{\f,\da\}$, is an ordered sequence, so that a labeling consists of a
bijection of sets $\{(x_1,1),..., (x_n,n) \}\to \{$incoming
edges$\}$, respecting the colors. We use the convention that, when
drawing a labeled tree in the plane, the top inputs from left to
right represent the sequence $\vec X$ of the labeling. Finally, we
say that a vertex is {\em binary} if it has valence $3$. A tree is
called {\em binary} if all of its vertices have valence $1$ or are
binary.
\end{Def}

\begin{Def}[Determinants]
We denote by $\mathrm{Edge}(T)$ the set of all edges of a tree $T$,
except the root, and by $|\mathrm{Edge}(T)|$ the cardinality of this
set. We define $\Det(T):= \bigwedge^{|\mathrm{Edge}(T)|}
k^{\mathrm{Edge}(T)}$. Notice that if $T$ is a tree with $n$ inputs
but no internal edges and if we let $S_n$ act on the tree by
permuting the inputs, then $\Det(T)$ becomes the sign representation
$sgn_n$ of $S_n$.
\end{Def}

\begin{Def}[Twisted dual]
Let $sgn_n$ be the sign representation of $S_n$. If the vector space
$V$ is an $S_n$-module, then the dual vector space $V^*$ is also an
$S_n$-module with the transposed action, and the space $V^\vee :=
V^*\otimes sgn_n$ is called the {\em twisted dual} or, using the
terminology of \cite{MSS}, the {\em Czech dual} of $V$.
\end{Def}

\begin{Def}[$\mathcal P(T)$]
Suppose $\mathcal P$ is a 2-colored operad, $T$ is 2-colored tree an
$v$ is a vertex of $T$. If $v$ has $n$ ``full'' inputs, $m$ ``dashed''
inputs and output color $x_v$, we set
\[
\mathcal P(v):= \left( \bigoplus_{\vec X\to
\mathrm{In}(v)}\mathcal P(\vec X;x_v) \right)_{S_{n+m}}
\]
where $(\cdot)_{S_{n+m}}$ denotes the coinvariants. The sum above is
taken over all $\vec X$ with $n$ many ``full'' colors and $m$ many
``dashed'' colors, and all labelings $\vec X\to In(v)$ respecting the
colors. Then we define
\[
\mathcal P(T):=\bigotimes_{ \text{vertex } v\text{ of  }T\text{ of
valence }\geq 3} \mathcal P(v).
\]
\end{Def}

\begin{Def}[Cobar Dual]
Given a 2-colored operad $\mathcal P$, we define its {\em cobar dual
2-colored operad} $\textbf{D} (\mathcal P)$ to be the differential
graded operad, given by the complexes $\textbf{D} (\mathcal
P)(\vec X;x):=$
\begin{multline}\label{DP}
\bigoplus_{
\substack{
  \text{trees  $T$ of}\\
  \text{type $(\vec X;x)$},\\
  \text{no internal edge}}}
\mathcal P(T)^*\otimes \Det(T)\stackrel{\de}{\lora}
\bigoplus_{
\substack{
  \text{trees } T \text{ of}\\
  \text{type }(\vec X;x),\\
  \text{1 internal edge}}}
\mathcal P(T)^*\otimes \Det(T)\stackrel{\de}{\lora}\\
\lora \dots\stackrel{\de}{\lora} \bigoplus_{
\substack{
  \text{trees } T \text{ of}\\
  \text{type }(\vec X;x),\\
  \text{binary tree}}}
\mathcal P(T)^*\otimes \Det(T).
\end{multline}
The above sum is always over isomorphism classes of colored, labeled
trees $T$. The operator $\de$ is given by the dual of the operad maps
\eqref{diopmaps}, which crunches an edge, twisted by a sign due to
the terms $\Det(T)$. Notice that this makes sense, since crunching an
edge preserves the type $(\vec X;x)$ of a tree $T$. The twisted sign
guarantees that $\de^2=0$. The degree is so that summands
corresponding to trees with $j$ internal edges have degree
$j+2-|\vec X|$. Thus the furthest right term has degree $0$. The
$S_n$ action is given by the $S_n$ action of the trees $T$ and the
operad maps are simply the attachments of trees.\\ Observe that since
there is only one tree with no internal edges, the furthest left term
in $\textbf{D} (\mathcal P)(\vec X;x)$ is just $\mathcal
P(\vec X;x)^\vee$.
\end{Def}

\subsection{Quadratic Koszul Colored Operads}

We want to define the notion of a quadratic colored operad, which by
definition is given by generators with 2 inputs and relations in 3
inputs. Then the quadratic dual will be defined by taking the
``orthogonal'' of the relations.\\
Let $E=\{E^{y,z}_x\}_{x,y,z\in\{\ff,\dda\}}$ be a set of $k$-vector
spaces, together with an $S_2$-action compatible with the colors.
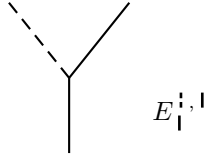
\[
\begin{pspicture}(0,0.5)(4,3.5)
 \psline[linestyle=dashed](2,2)(1.2,3)
 \psline(2,2)(2.8,3)
 \psline(2,2)(2,1)
 \rput[b](3.5,1.3){$E^{\dda,\ff}_{\ff}$}
\end{pspicture}
\]
We want $E$ to be the binary generating set of a 2-colored operad,
where $x$, $y$ and $z$ correspond to the colors of a binary
vertex. Let $v$ be a binary vertex of type $(y,z;x)$ in
a binary tree $T$. Then, the expression $E(v)$ is defined to be
\[
E(v):= \left( E^{y,z}_x \oplus E^{z,y}_x \right)_{S_2}
\]
and we set $E(T):=\bigotimes_{\text{binary vertex }v\text{ of }T} E(v)$.
With this notation, we define the {\em free 2-colored operad generated by}
$E$ to be
\[
\mathcal F(E)(\vec X;x):= \bigoplus_{
  \text{binary trees $T$ of type }(\vec X;x)
} E(T).
\]
The $S_n$-action is given by an obvious permutation of the tree using
the $S_2$-action on $E$, and the composition maps are given by
attaching trees. This definition can readily be seen to define a
2-colored operad.\\
An {\em ideal} $\mathcal I$ of a 2-colored operad $\mathcal P$ is a
collection of $S_n$-sub-modules $\mathcal I(\vec X;x)\subset
\mathcal P(\vec X;x)$ such that $f\circ_{i} g$ belongs to the ideal
whenever $f$ or $g$ or both belong to the ideal.

\begin{Def}[Quadratic Operad]
A 2-colored operad $\mathcal P$ is said to be {\em quadratic} if
$\mathcal P=\mathcal F(E)/(R)$ where $\mathcal F(E)$ is the free
2-colored operad on some generators $E$, and $(R)$ is the ideal in
$\mathcal F(E)$ generated by a subspace with 3 inputs, called the
{\em relations}: $$R\subset \bigoplus_{|\vec X|=3,\ x=\{\ff,\dda\}}
\mathcal F(E)(\vec X;x).$$
\end{Def}

Recall that if the vector space $V$ is an $S_n$-module and $sgn_n$ is the
sign representation, then we defined $V^\vee$ to be $V^*\otimes sgn_n$.

\begin{Def}[Quadratic Dual]
For every quadratic 2-colored operad $\mathcal P$, we define
the {\em quadratic dual 2-colored operad} $\mathcal P^!:=\mathcal
F(E)^\vee/(R^\perp)$, where $(R^\perp)$ is the ideal in $\mathcal
F(E)^\vee$ generated by the orthogonal complement $R^\perp$ of $R$
in $$\left(\bigoplus_{|\vec X|=3,\ x\in\{\ff,\dda\}} \mathcal
F(E)(\vec X;x) \right) ^\vee.$$
\end{Def}

Notice that $\mathcal F(E)(\vec X;x)^\vee = \mathcal F(E^\vee) (\vec
X;x)$ (see \cite{GK}), so that $\mathcal P^!$ is generated by
$E^\vee$ with relations $R^\perp$. Let us describe the composition
maps of $\mathcal P^!$, i.e., the maps of equation \eqref{diopmaps}.
They are given by the composition of the following maps:
\begin{equation}
\mathcal F(E^\vee)(\vec X;x) \otimes \mathcal F(E^\vee)(\vec Y;y) \to
\mathcal F(E^\vee)(\vec Z;x) \to \mathcal F (E^\vee)/(R^\perp)(\vec
Z;x) = \mathcal P^!(\vec Z;x),
\end{equation}
where the first map is the tensor product made according to the trees
in $\mathcal F(E^\vee)$, and the second is the natural projection.
Then, one can check that these maps descend to the quotient $\mathcal
F(E)^\vee/(R^\perp)(\vec X;x) \otimes \mathcal F(E)^\vee/
(R^\perp)(\vec Y;y) = \mathcal P^!(\vec X;x) \otimes \mathcal
P^!(\vec Y;y)$.

\begin{Def}[Koszulness]
If $\mathcal P=\mathcal F(E)/(R)$ is a quadratic 2-colored operad,
then we can consider the cobar dual 2-colored operad $\textbf{D}
(\mathcal P^!)$ of its quadratic dual $\mathcal P^!$. In particular
the term of degree zero in the complex $\textbf{D} (\mathcal
P^!)(\vec X;x)$ is given by
\[
\textbf{D}(\mathcal P^!)(\vec X;x)^0 = \bigoplus_{
\substack{\text{binary trees} \\
\text{of type }(\vec X;x)}}
E(T) =\mathcal F(E)(\vec X;x)
\]
and therefore we have a sequence of maps
\begin{multline}\label{sequence_H0}
\cdots \stackrel{\de}{\lora} \textbf{D}(\mathcal P^!)
(\vec X;x)^{-1} \stackrel{\de}{\lora} \textbf{D}(\mathcal P^!)
(\vec X;x)^0=\\=\mathcal F(E)(\vec X;x) \stackrel{proj}\lora
\mathcal F(E)/(R)(\vec X;x)=\mathcal P(\vec X;x)
\end{multline}
Notice that all the vertices of trees in the term
$\textbf{D}(\mathcal P^!)(\vec X;x)^{-1}$ are binary except one
ternary vertex, and that the space $\mathcal P^!(\vec X;x)$ is
the dual of the space $R$ of relations of $\mathcal P$. Thus the
image of $\textbf{D}(\mathcal P^!)(\vec X;x)^{-1}$ in
$\textbf{D}(\mathcal P^!)(\vec X;x)^{0}$ is the space (of
consequences) of the relations in $\mathcal P(\vec X;x)$. This
means that the morphism $proj$ induces an isomorphism
\begin{equation}\label{koszul_at_0}
H^0(\textbf{D}(\mathcal P^!)(\vec X;x))\cong \mathcal
P(\vec X;x).
\end{equation}
A quadratic operad is said to be {\em Koszul} if the cobar dual
complex $\textbf{D}(\mathcal P^!)(\vec X;x)$ is quasi-isomorphic to
$\mathcal P(\vec X;x)$, i.e., by \eqref{koszul_at_0}, has homology
concentrated in degree zero, for any $\vec X$ and $x$.
\end{Def}

\begin{Rem} Another excellent source for Koszul duality of colored
operads, which includes all the expected standard results of the
theory, can be found in Pepijn van der Laan's Ph.D. Thesis \cite{L1}
(see also \cite{L2}).
\end{Rem}

\subsection{Algebras and Infinity Algebras, Modules and Infinity Modules}
\label{ssec:am}

We now want to describe an object which encodes the structure of
modules over an operad $\mathcal O$. We have to use the language of
2-colored operads as follows

\begin{Def}[$\overline{\mathcal O}$]\label{def_O_bar}
Suppose $\mathcal O$ is a 1-colored operad, then define
\[
 \overline{\mathcal O}(\vec X;x):=
  \begin{cases}
    \mathcal O(n) & \text{if $x$ is ``full'', and $\vec X=(\f,\ldots,\f)$,} \\
    \mathcal O(n) & \text{if $x$ is ``dashed'', and $\vec X$ has exactly
    one ``dashed'' input,} \\
    \{0\} & \text{otherwise}.
  \end{cases}
\]
The $S_n$-action on $\overline{\mathcal O}$ is defined to be the
$S_n$-action on $\mathcal O$ together with a possible change of the
coloring according to \eqref{equimap}. Similarly, the composition maps
are given by the composition maps from $\mathcal O$, again with
additionally keeping track of the coloring.
\end{Def}

Thus, the only non-empty spaces are those which consist only of the
``full'' color (and thus including the operad $\mathcal O$ itself)
and those which have exactly one incoming ``dashed'' color together
with outgoing ``dashed'' color. The idea behind the definition is that
the ``full''-colored multiplications will describe algebras over
$\mathcal O$, whereas the ``dashed'' part will describe modules over
algebras over $\mathcal O$.\\
If $\mathcal O$ happens to be quadratic and Koszul, then so does
$\overline{\mathcal O}$:

\begin{Lem}\label{O_bar_quadratic}
Let $\mathcal O$ be quadratic with generators $E=E^{\ff,\ff}_{\ff}$ and
relations $R\subset\mathcal F(E)(3)$. Then $\overline{\mathcal O}$
is quadratic with generators $\overline{E}= \{\overline{E}^{y,z}_x\}$,
given by $\overline{E}^{\ff,\ff}_{\ff}:=E^{\ff,\ff}_{\ff}$,
$\overline{E}^{\ff,\dda}_{\dda}:=E^{\ff,\ff}_{\ff}$,
$\overline{E}^{\dda,\ff}_{\dda}:=E^{\ff,\ff}_{\ff}$, and $\{0\}$ otherwise.
The relations are given by $R$, whenever the space $\mathcal
F(\overline{E})$ is non-empty:
\begin{eqnarray*}
R\subset\mathcal F(E)(3)\cong\mathcal F(\overline{E})(\f,\f,\f;\f), \\
R\subset\mathcal F(E)(3)\cong\mathcal F(\overline{E})(\f,\f,\da;\da), \\
R\subset\mathcal F(E)(3)\cong\mathcal F(\overline{E})(\f,\da,\f;\da), \\
R\subset\mathcal F(E)(3)\cong\mathcal F(\overline{E})(\da,\f,\f;\da).
\end{eqnarray*}
\end{Lem}

\begin{proof} Notice that the free operad generated by
$\overline{E}$ is simply $\mathcal F(E)$ together with copies of
itself for every possible input position of one ``dashed'' input
which then necessarily has a ``dashed'' output. For each of those
copies, we take the same relations from $\mathcal O$, effectively
giving $n+1$ copies of $\mathcal O$ for the $n$-input maps of
$\overline{\mathcal O}$.
\end{proof}
\begin{Prop}\label{O_bar_Koszul}
If $\mathcal O$ is a quadratic and Koszul, then $\overline{\mathcal
O}$ has a resolution given by
\[
\textbf{D}(\overline{\mathcal O^!}) (\vec X;x)\to \overline{
\mathcal O}(\vec X;x)=\mathcal O(n)
\]
\end{Prop}
\begin{proof} The proof uses the same idea as the proof of Lemma
\ref{O_bar_quadratic}, namely the arguments for $\mathcal O$ split
into several copies with either only ``full'' colors or one
``dashed'' input color. In other words, if we ignore the colors of
the map $\textbf{D}(\overline{\mathcal O^!})(\vec X;x) \to
\overline{\mathcal O}(\vec X;x),$ we get exactly $\textbf{D}(\mathcal
O^!)(n)\to \mathcal O(n)$, where $|\vec X|=n$.
\end{proof}
The above proposition can also be found in Pepijn van der Laan's
Ph.D. Thesis \cite{L1} Theorem 3.4.2, and in \cite{L2}.

We now consider algebras over the 2-colored operad $\overline{\mathcal
O}$. First recall that an algebra over $\mathcal O$ consists of maps
$\mathcal O(n)\to Hom(A^{\otimes n},A)$, i.e., $\mathcal O(n) \otimes
A^{\otimes n}\to A$, which are associative and $S_n$-equivariant.
Similarly, when applying Definition \eqref{def_col_alg} to
$\overline{\mathcal O}$, we see that algebras over $\overline{\mathcal
O}$ are given by a pair of $k$-vector spaces $A$ and $M$ together with
maps $\mathcal O(n) \otimes A^{\otimes n}\to A$ and
$\bigoplus_{r+s=n-1} \mathcal O(n) \otimes A^{\otimes r}\otimes
M\otimes A^{\otimes s}\to M$ being again associative and
$S_n$-equivariant. We will sometimes abuse notation and call $A$
an {\em algebra over} $\mathcal O$ and $M$ a {\em module over} $A$
{\em over} $\mathcal O$.\\
Our aim in the rest of the subsection is to identify algebras over
$\textbf{D}(\mathcal O)$ and $\textbf{D}(\overline{\mathcal O})$ in
the case when $\mathcal O$ is a quadratic operad generated by
$E=\mathcal O (2)$ and $\overline{\mathcal O}$ is the quadratic
2-colored operad generated by $\overline{E}$, defined above.

\begin{Def}[Free Algebra]\label{freealgebra}
For an operad $\mathcal O$, the free $\mathcal O$-algebra generated
by the vector space $A$ is defined as $$F_{\mathcal O} A :=
\bigoplus_{n\geq 1} \left(\mathcal O(n)\otimes A^{\otimes
n}\right)_{S_n}.$$ $F_{\mathcal O} A$ is an algebra over $\mathcal
O$, i.e., there are maps $ \gamma:\mathcal O(k)\otimes (F_{\mathcal O}
A )^{\otimes k} \to F_{\mathcal O} A$ coming from the composition in
$\mathcal O$ and the tensor products, which satisfy the required
compatibility conditions.
\end{Def}

\begin{Def}[Free Module]\label{freemodule}
If $A$ and $M$ are $k$-vector spaces, we define
\[
F_{\mathcal O,A} M := \bigoplus_{n\geq 1}
\left(\bigoplus_{r+s=n-1} \mathcal O(n) \otimes A^{\otimes r}
\otimes M \otimes A^{\otimes s}\right)_{S_{n}}.
\]
Then $F_{\mathcal O,A} M$ is a module over $F_{\mathcal O}A$ over
$\mathcal O$, which means that there are maps $\gamma^M:\bigoplus_{r+s=n-1}
\mathcal O(n)\otimes (F_{\mathcal O} A )^{\otimes r} \otimes
F_{\mathcal O,A} M \otimes (F_{\mathcal O} A )^{\otimes s} \to
F_{\mathcal O,A} M$ given by composition of elements of the operad
and tensor product of elements of $A$. Again, these maps satisfy the
required axioms (see \cite{GK}).\\ Notice that the coinvariants
$(\bigoplus_{r+s=n-1} \mathcal O(n) \otimes A^{\otimes r} \otimes M
\otimes A^{\otimes s})_{S_{n}}$ are canonically isomorphic to the
coinvariants $(\mathcal O(n) \otimes A^{\otimes n-1} \otimes
M)_{S_{n-1} \times id}$, where $S_{n-1} \times id$ is the subgroup of
$S_n$ fixing the last element. We will switch between these two
descriptions of $F_{ \mathcal O, A}M$, according to whichever one
seems more useful.
\end{Def}

\begin{Def}[Algebra Derivation]
We define a {\em derivation} $d\in \mathrm{Der}(F_{\mathcal O}A)$ to
be a map from $F_{\mathcal O} A$ to itself, making the following
diagram commute:
\[
\begin{CD}
\mathcal O(k) \otimes (F_{\mathcal O}A)^{\otimes k}& @>\gamma>> &
F_{\mathcal O}A\\ @V\sum_i \mathrm{id}\otimes \mathrm{id}^{\otimes i}
\otimes d \otimes \mathrm{id}^{\otimes (k-i-1)}VV & & @VVdV\\
\mathcal O(k) \otimes (F_{\mathcal O}A)^{\otimes k}& @>\gamma>> &
F_{\mathcal O}A
\end{CD}
\]
where $\gamma$ is the $\mathcal O$-algebra structure of $F_{\mathcal
O} A$.\\ Notice that, if $\mathcal O$ is quadratic, then this
condition is equivalent to saying that for every $\mu\in\mathcal O
(2)$ and $a, b\in F_{\mathcal O}A$ one has $
d(\mu(a,b))=\mu(d(a),b)+(-1)^{|a|\cdot|d|}\mu(a,d(b))$, where by the
symbol $\mu(a,b)$ we mean $\gamma(\mu, a, b)$.\\ Also notice that a
derivation $d$ is completely determined by its restriction to $A\cong
\mathcal O(1)\otimes A \subset F_{\mathcal O} A$, i.e., by maps
$$d_n:A \to (\mathcal O(n) \otimes A^{\otimes n})_{S_n},$$ for
$n\in\N$.
\end{Def}

\begin{Def}[Module Derivation]
Given a derivation $d\in \mathrm{Der}(F_{\mathcal O}A)$, we define a
{\em derivation} $g\in \mathrm{Der}_d (F_{\mathcal O, A}M)$ {\em
over} $d$ to be a map from $F_{\mathcal O, A}M$ to itself, making the
following diagram commutative:
\[
\begin{CD}
\mathcal O(k) \otimes (F_{\mathcal O}A)^{\otimes (k-1)}\otimes
F_{\mathcal O, A}M & @>\gamma^M>> & F_{\mathcal O, A}M\\ @V\alpha VV
& & @VVgV\\ \mathcal O(k) \otimes (F_{\mathcal O}A)^{\otimes
(k-1)}\otimes F_{\mathcal O, A}M & @>\gamma^M>> & F_{\mathcal O, A}M
\end{CD}
\]
Here $\alpha=\sum_{i<k} \mathrm{id}\otimes \mathrm{id}^{\otimes i}
\otimes d \otimes \mathrm{id}^{\otimes (k-i-1)}+ \mathrm{id}\otimes
\mathrm{id}^{\otimes (k-1)} \otimes g$, and $\gamma^M$ is the module
structure of $F_{\mathcal O, A} M$ over $F_{\mathcal O}A$ over
$\mathcal O$.\\ Notice that, if $\mathcal O$ is quadratic, then this
condition is equivalent to saying that for every $\mu\in\mathcal O
(2)$, $a\in F_{\mathcal O}A$ and $m\in F_{\mathcal O, A}M$, one has $
g(\mu(a,m))=\mu(d(a),m)+ (-1)^{|a|\cdot|g|}\mu(a,g(m))$, where by the
symbol $\mu(a,m)$ we mean $\gamma^M(\mu, a, m)$.\\ Again, a
derivation $g$ over $d$ is completely determined by its restriction
to $M\cong \mathcal O(1)\otimes M \subset F_{\mathcal O, A} M$, i.e.,
by maps $$g_{n+2}: M \to \left(\mathcal O(n+1) \otimes (A^{\otimes
n}\otimes M)\right)_{S_{n}\times id},$$ for $n\in\N$, and where by
$S_{n}\times id$ we mean the subgroup of $S_{n+1}$ whose elements
leave the last component fixed.
\end{Def}

With this we can identify algebras over $\textbf{D} (\mathcal O)$ and
$\textbf{D} (\overline{\mathcal O})$. In order to state the Theorem
we need the following
\begin{Def}[Suspension]
Let $A=\oplus_{j\in \Z} A_j$ be a graded vector space over $k$. We
define the {\em suspension}\/ $sA$ of $A$ as the graded vector space
$sA=\oplus_{j\in \Z} (sA)_j$ with $(sA)_j=A_{j-1}$. Therefore the
suspension map $s\colon A\to sA$ sending $v$ to $sv:=v$ is an
isomorphism of degree +1. Similarly, we set $s^{-1}A = \oplus_{j\in \Z}
(s^{-1}A)_j$ with $(s^{-1}A)_j=A_{j+1}$.
\end{Def}

\begin{Thm}[\cite{GK}, Proposition (4.2.15)]\label{regular_koszuality}
Let $(A,d_A)$ be a differential graded space. Then, giving $A$ the
structure of an algebra over $\textbf{D} (\mathcal O)$ is the same as
giving a derivation $d\in \mathrm{Der} (F _{\mathcal O}\,sA^*)$ of
degree $1$, with $d^2=0$.
\end{Thm}
If $\mathcal O$ is quadratic Koszul, then by definition it is
$\textbf{D}(\mathcal O^!)\cong \mathcal O$, and with this Theorem, we
define a {\em homotopy algebra over} $\mathcal O$, or {\em homotopy
$\mathcal O$-algebra}, to be given by a differential $d\in
\mathrm{Der}(F _{\mathcal O^!}\,sA^*)$ with $d^2=0$.

\begin{Thm}\label{O_bar_algebras} Let $(A,d_A)$ and $(M,d_M)$ be
differential graded spaces. Then giving $(A,M)$ the structure of an
algebra over $\textbf{D} (\overline{\mathcal O})$ is the same as
giving a differential $d\in \mathrm{Der}(F _{\mathcal O}\,sA^*)$ of
degree $1$, with $d^2=0$, together with a differential $g\in
\mathrm{Der}_d (F_{\mathcal O, sA^*}sM^*)$ over $d$ of degree $1$,
with $g^2=0$.
\end{Thm}
If $\mathcal O$ is quadratic Koszul, then by Proposition
\ref{O_bar_Koszul} it is $\textbf{D}(\overline{\mathcal O^!})\cong
\overline{\mathcal O}$, and with this Theorem, we define a {\em
homotopy module over the homotopy algebra} $A$ {\em over} $\mathcal
O$, or {\em homotopy $\mathcal O$-module}, to be given by
differentials $d\in \mathrm{Der}(F _{\mathcal O^!}\,sA^*)$ and $g\in
\mathrm{Der}_d (F_{\mathcal O^!, sA^*}sM^*)$ over $d$, with $d^2=0$
and $g^2=0$.\\

The proof for Theorem~\ref{regular_koszuality} can be found in
\cite{GK}. The proof for Theorem~\ref{O_bar_algebras} is analogous to
the one of Theorem \ref{regular_koszuality}, with the difference
that we now have to keep track of the ``full'' and the ``dashed''
colors, similar to Proposition \ref{O_bar_Koszul}. We will sketch
part of these proofs when dealing with the analogous statement for
homotopy inner products in Theorem~\fullref{O_hat_algebras}.

\section{Inner Products over Cyclic Operads}
\label{sec:cy}

We first recall the definition of cyclic operads $\mathcal O$ and
associate to it an object $\widehat{\mathcal O}$ which describes
invariant inner products over $\mathcal O$. Then, in subsections
\ref{quadrat-koszul} and \ref{proof_main}, we show how quadraticity
and Koszulness of $\mathcal O$ implies quadraticity and Koszulness of
$\widehat{\mathcal O}$. Finally in subsections \ref{homotop-ip} and
\ref{Ex:A,C,L}, we identify homotopy algebras over $\widehat{\mathcal
O}$ in terms of derivations and module maps.

\subsection{Cyclic Operads and $\widehat{\mathcal O}$}
\label{cyclic-op}

\begin{Def}[Cyclic Operad]
Let $\mathcal O$ be an operad. Recall from \cite{GeK} and \cite{MSS}
that $\mathcal O$ is called {\em cyclic}, if there is an action of
the symmetric group $S_{n+1}$ on $\mathcal O(n)$, which extends the
given $S_n$-action, and satisfies, for $1\in\mathcal O(1)$,
$\alpha\in \mathcal O(m)$, $\beta\in \mathcal O(n)$ the following
relations:
\begin{eqnarray}
\label{compos_cyclic1} \tau_2(1)&=&1,\\ \label{compos_cyclic2}
\tau_{m+n}(\alpha\circ_k \beta)&=&\tau_{m+1}(\alpha)\circ_{k+1}
\beta,\quad\quad\quad \text{ for } k<m \\ \label{compos_cyclic3}
\tau_{m+n}(\alpha\circ_m \beta)&=&\tau_{n+1}(\beta)\circ_1
\tau_{m+1}(\alpha),
\end{eqnarray}
where $\tau_{j}\in S_{j}$ denotes the cyclic rotation of $j$ elements
$\tau_{j} :=1\in\Z_{j}\subset S_{j}$.
\end{Def}

We want to use this extra datum to define a new object $\widehat{
\mathcal O}$ associated to $\mathcal O$. Informally, $\widehat{
\mathcal O}$ uses the $S_{n+1}$ action to turn the output of a tree
into another input, and we thus have to look trees with inputs but no
output.
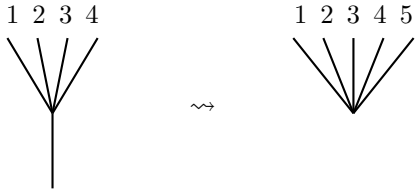
\[
\begin{pspicture}(0,0)(4,4)
 \psline(2,2)(1.4,3)
 \psline(2,2)(1.8,3)
 \psline(2,2)(2.2,3)
 \psline(2,2)(2.6,3)
 \psline(2,2)(2,1)
 \rput[b](2,3.2){$1\,\,\, 2\,\,\, 3\,\,\, 4$}
 \rput[b](4,2){$\rightsquigarrow$}
\end{pspicture}
\begin{pspicture}(0,0)(4,4)
 \psline(2,2)(1.2,3)
 \psline(2,2)(1.6,3)
 \psline(2,2)(2,3)
 \psline(2,2)(2.4,3)
 \psline(2,2)(2.8,3)
 \rput[b](2,3.2){$1\,\,\, 2\,\,\, 3\,\,\, 4\,\,\, 5$}
\end{pspicture}
\]
Therefore, we need the following generalization of 2-colored operads.

\begin{Def}[0/1-Operad]\label{0/1-operad}
A {\em 0/1-operad}\/ is a sub-operad of a 3-colored operad with colors
$\{\f,\da,\varnothing \}$ such that the color $\varnothing$ can appear
only as an output. As in the 2-colored case, we require that only
nontrivial spaces with one input are $\mathcal P(\f;\f)=k$ and
$\mathcal P(\da;\da)=k$, and we fix generators of these spaces.
\end{Def}

In practice, a 0/1-operad consists of $k$-vector spaces
$\mathcal P(\vec X ;x)$, where $\vec X$ indicates, as before, which
inputs are colored with the ``full'' or the ``dashed'' color and
$x\in\{\f,\da,\varnothing\}$ denotes the output color. Since the
color $\varnothing$ cannot appear as an input, we may use the
following convention: we represent the $\varnothing$ output with a
blank line, i.e., with no line, and we say that the operation
``has no output''.
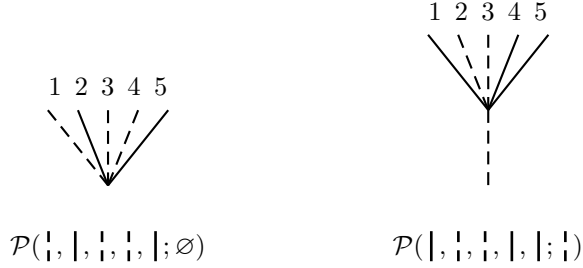
\[
\begin{pspicture}(0,1)(4,4)
 \psline[linestyle=dashed, arrowsize=0.1, arrowinset=0](2,2)(1.2,3)
 \psline[arrowsize=0.1, arrowinset=0](2,2)(1.6,3)
 \psline[linestyle=dashed, arrowsize=0.1, arrowinset=0](2,2)(2,3)
 \psline[linestyle=dashed, arrowsize=0.1, arrowinset=0](2,2)(2.4,3)
 \psline[arrowsize=0.1, arrowinset=0](2,2)(2.8,3)
 \rput[b](2,1){$\mathcal P (\da,\f,\da,\da,\f;\varnothing)$}
 \rput[b](2,3.2){$1\,\,\, 2\,\,\, 3\,\,\, 4\,\,\, 5$}
\end{pspicture}
\quad \quad \quad
\begin{pspicture}(0,0)(4,3.6)
 \psline[arrowsize=0.1, arrowinset=0](2,2)(1.2,3)
 \psline[linestyle=dashed, arrowsize=0.1, arrowinset=0](2,2)(1.6,3)
 \psline[linestyle=dashed, arrowsize=0.1, arrowinset=0](2,2)(2,3)
 \psline[arrowsize=0.1, arrowinset=0](2,2)(2.4,3)
 \psline[arrowsize=0.1, arrowinset=0](2,2)(2.8,3)
 \psline[linestyle=dashed, arrowsize=0.1, arrowinset=0](2,2)(2,1)
 \rput[b](2,0){$\mathcal P (\f,\da,\da,\f,\f;\da)$}
 \rput[b](2,3.2){$1\,\,\, 2\,\,\, 3\,\,\, 4\,\,\, 5$}
\end{pspicture}
\]

Clearly we can still form the cobar dual $\textbf{D}(\mathcal P)$.
More explicitly, $\textbf{D}(\mathcal P)(\vec X;\f)$ and
$\textbf{D}(\mathcal P)(\vec X;\da)$ are the complexes of
equation~\eqref{DP}, and similarly $\textbf{D} (\mathcal
P)(\vec X;\varnothing)$ is the complex of equation~\eqref{DP} where
we change the trees to be of type $(\vec X;\varnothing)$, i.e., trees
with ``full'' and ``dashed'' edges but with no output.\\ The
canonical example of this structure is the endomorphism 0/1-operad
given for $k$-vector spaces $A$ and $M$ by
\begin{align*}
\End^{A,M}(\vec X;\f)=& Hom(\text{tensor products of $A$ and $M$},A)\\
\End^{A,M}(\vec X;\da)=& Hom(\text{tensor products of $A$ and $M$},M)\\
\End^{A,M}(\vec X;\varnothing)=& Hom(\text{tensor products of $A$ and
$M$},k).
\end{align*}
By a slight abuse of language we will say that the pair $(A,M)$ is an
algebra over the 0/1-operad $\mathcal P$ if there exists a 0/1-operad
map  $\mathcal P \to \End^{A,M}$.\\
Now, we will use this new concept to define, for a given cyclic
operad $\mathcal O$, the 0/1-operad $\widehat{\mathcal O}$.

\begin{Def}[$\widehat{\mathcal O}$]
\label{def_O_hat}
Let $\mathcal O$ be a cyclic operad. If $|\vec X|=n$, then let
\[
 \widehat{\mathcal O}(\vec X;x):=
   \begin{cases}
    \mathcal O(n) & \text{if } x \text{ is ``full'', and } \vec X=(\f,\ldots,\f),\\
    \mathcal O(n) & \text{if } x \text{ is ``dashed'', and $\vec X$ has
    exactly one ``dashed'' input}\\
    \mathcal O(n-1) & \text{if } x=\varnothing \text{ and
    $\vec X$ has exactly two ``dashed' inputs,} \\
    \{0\} & \text{otherwise}.
  \end{cases}
\]
The definition of $\widehat{\mathcal O}(\vec X,\varnothing)$ is
motivated by the idea that one considers trees with $n-1$ inputs and
one output, and then turns this output into a new input:
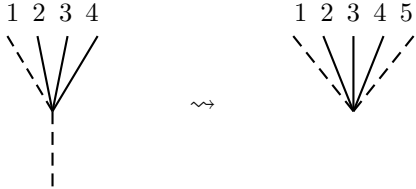
\[
\begin{pspicture}(0,0)(4,4)
 \psline[linestyle=dashed](2,2)(1.4,3)
 \psline(2,2)(1.8,3)
 \psline(2,2)(2.2,3)
 \psline(2,2)(2.6,3)
 \psline[linestyle=dashed](2,2)(2,1)
 \rput[b](2,3.2){$1\,\,\, 2\,\,\, 3\,\,\, 4$}
 \rput[b](4,2){$\rightsquigarrow$}
\end{pspicture}
\begin{pspicture}(0,0)(4,4)
 \psline[linestyle=dashed](2,2)(1.2,3)
 \psline(2,2)(1.6,3)
 \psline(2,2)(2,3)
 \psline(2,2)(2.4,3)
 \psline[linestyle=dashed](2,2)(2.8,3)
 \rput[b](2,3.2){$1\,\,\, 2\,\,\, 3\,\,\, 4\,\,\, 5$}
\end{pspicture}
\]
We define the $S_n$-action on $\widehat{\mathcal O}(\vec X;\f)$ and
$\widehat{\mathcal O}(\vec X;\da)$ as before by using the
$S_n$-action on $\mathcal O(n)$, and the $S_n$-action on
$\widehat{\mathcal O} (\vec X;\varnothing)$ by using the $S_n$-action
on $\mathcal O(n-1)$ given by the cyclicity of $\mathcal O$.\\
Diagrams with different positions of the two ``dashed'' inputs can be
mapped to each other using the action of the symmetric group. In
fact, as each $\sigma\in S_{n+1}$ induces an isomorphism which
preserves all the structure, any statement about diagrams with a
fixed choice of position of ``dashed'' inputs immediately carries
over to any other choice of positions of ``dashed'' inputs. We
therefore often restrict our attention to the choice where the two
``dashed'' inputs are at the far left and the far right, as shown in
the above picture.\\ It is left to define the composition. On
$\widehat{\mathcal O}(\vec X;\f)$ and $\widehat{\mathcal O}(\vec
X;\da)$, the composition is simply the composition in $\mathcal
O(n)$, so that it clearly satisfies equivariance and associativity.
If $|\vec X|=n+1$, then on $\widehat{\mathcal O}(\vec
X;\varnothing)=\mathcal O(n)$, the composition is predetermined on
the first $n$ components by the usual composition in $\mathcal O$. As
for the last component, we define
\begin{equation}\label{def_cyclic_compos}
\alpha\circ_{m+1} \beta:=\tau_{n+m} (\tau^{-1}_{m+1}(\alpha)\circ_m
\beta)\stackrel{\eqref{compos_cyclic3}}{=} \tau_{n+1}(\beta) \circ_1
\alpha
\end{equation}
\[
\begin{pspicture}(1,1.6)(10.5,5.6)
 \psline[linestyle=dashed](2,2)(1.2,2.9)
 \psline(2,2)(1.6,2.9)
 \psline(2,2)(2,2.9)
 \psline(2,2)(2.4,2.9)
 \psline[linestyle=dashed](2,2)(2.8,2.9)
 \rput[b](3,2){$\alpha$} \rput[b](2.4,3.4){$\beta$}
 \psline[linestyle=dashed](2.8,3)(2.8,3.5)
 \psline[linestyle=dashed](2.8,3.5)(3,4)
 \psline(2.8,3.5)(2.8,4)
 \psline(2.8,3.5)(2.6,4)
 \rput[b](3.5,3){$:=$}
 \psline[linestyle=dashed](5,2)(4.2,2.9)
 \psline[linestyle=dashed](4.2,3)(4.3,3.1)(5,3.25)(5.7,3.4)(5.8,3.5)
 \psline(5,2)(4.6,2.9)                   \psline(4.6,3)(4.2,3.5)
 \psline(5,2)(5,2.9)                     \psline(5,3)(4.6,3.5)
 \psline(5,2)(5.4,2.9)                   \psline(5.4,3)(5,3.5)
 \psline[linestyle=dashed](5,2)(5.8,2.9) \psline[linestyle=dashed](5.8,3)(5.4,3.5)
 \rput[b](6,2){$\alpha$} \rput[b](6.4,3.1){$\tau_{m+1}^{-1}$}
 \psline[linestyle=dashed](5.4,3.6)(5.4,4.1)
 \psline[linestyle=dashed](5.4,4.1)(5.6,4.5)
 \psline(5.4,4.1)(5.4,4.5)
 \psline(5.4,4.1)(5.2,4.5)
 \rput[b](5.2,3.8){$\beta$}
 \psline[linestyle=dashed](5.8,3.6)(5.8,4.5)
 \psline(4.2,3.6)(4.2,4.5)
 \psline(4.6,3.6)(4.6,4.5)
 \psline(5,3.6)(5,4.5)
 \psline(4.2,4.6)(4.4,5.3)
 \psline(4.6,4.6)(4.8,5.3)
 \psline(5  ,4.6)(5.2,5.3)
 \psline(5.2,4.6)(5.4,5.3)
 \psline(5.4,4.6)(5.6,5.3)
 \psline[linestyle=dashed](5.6,4.6)(5.8,5.3)
 \psline[linestyle=dashed](5.8,4.6)(5.6,4.8)(5,4.95)(4.4,5.1)(4.2,5.3)
 \rput[b](6.4,4.7){$\tau_{n+m}$}
 \rput[b](7.3,3){$=$}
 \psline[linestyle=dashed](8.55,2)(8.1,2.9)
 \psline(8.55,2)(8.4,2.9)
 \psline(8.55,2)(8.7,2.9)
 \psline[linestyle=dashed](8.55,2)(9,2.9)
 \rput[b](9.7,2){$\tau_{n+1}(\beta)$}
 \psline[linestyle=dashed](8.1,3)(8.1,3.5)
 \psline[linestyle=dashed](8.1,3.5)(7.8,4)
 \psline(8.1,3.5)(8  ,4)
 \psline(8.1,3.5)(8.2,4)
 \psline(8.1,3.5)(8.4,4)
 \rput[b](8.5,3.5){$\alpha$}
\end{pspicture}
\]
where $\alpha\in \widehat{\mathcal O}(\vec X;\varnothing) =\mathcal
O(m)$ has $m+1$ inputs, and $\beta\in\widehat{\mathcal
O}(\vec Y;\da)=\mathcal O(n)$ (or similarly $\beta\in\widehat{\mathcal
O}(\vec Y;\f)$) has $n$ inputs. It is clear that this will satisfy
equivariance, since equivariance was just used to define the
composition. The next Lemma establishes the final property for
$\widehat{\mathcal O}$ being a 0/1-operad.
\begin{Lem}
The composition in $\widehat{\mathcal O}$ satisfies the
associativity axiom.
\end{Lem}
\begin{proof} By definition we have that the composition is
just the usual composition in $\mathcal O(n)$, except for inserting
trees in the last input of elements in $\widehat{\mathcal
O}(\vec X;\varnothing)$. Thus, except for composition in the last
spot, associativity of $\widehat{\mathcal O}$ follows from the
associativity of $\mathcal O$.\\%
Now, let $\alpha\in \widehat{\mathcal O}(\vec X;\varnothing)\cong
\mathcal O(m)$, $\beta\in \widehat{\mathcal O}(\vec Y;y)\cong
\mathcal O(n)$, and $\gamma\in \widehat{\mathcal O} (\vec Z;z) \cong
\mathcal O(p)$, where $y,z\in\{\f,\da\}$. Then, associativity is
satisfied, because for $1\leq j\leq m$, it is
\begin{multline*}
(\alpha\circ_{m+1}\beta)\circ_{j} \gamma
\stackrel{\mathit{\eqref{def_cyclic_compos}}}{=} (\tau_{n+1}(\beta)
\circ_1 \alpha)\circ_j \gamma \stackrel{\mathit{op.comp}}{=}\\
=\tau_{n+1}(\beta) \circ_1 (\alpha\circ_j \gamma )
\stackrel{\mathit{\eqref{def_cyclic_compos}}}{=} (\alpha\circ_j
\gamma )\circ_{m+p}\beta,
\end{multline*}
and for $m< j< m+n$, it is
\begin{multline*}
(\alpha\circ_{m+1}\beta)\circ_{j} \gamma
\stackrel{\mathit{\eqref{def_cyclic_compos}}}{=} (\tau_{n+1}(\beta)
\circ_1 \alpha)\circ_j \gamma  \stackrel{\mathit{op.comp}}{=}
(\tau_{n+1}(\beta) \circ_{j-m+1} \gamma)\circ_1 \alpha
\stackrel{\eqref{compos_cyclic2}}{=}\\ = \tau_{n+p}(\beta
\circ_{j-m}\gamma)\circ_1\alpha \stackrel{\mathit{
\eqref{def_cyclic_compos}}}{=} \alpha\circ_{m+1}
(\beta\circ_{j-m}\gamma),
\end{multline*}
while
\begin{multline*}
(\alpha\circ_{m+1}\beta)\circ_{m+n} \gamma
\stackrel{\mathit{\eqref{def_cyclic_compos}}}{=} (\tau_{n+1}(\beta)
\circ_1 \alpha)\circ_{m+n} \gamma  \stackrel{\mathit{
\eqref{def_cyclic_compos}}}{=}\\ = \tau_{p+1}(\gamma)\circ_1
(\tau_{n+1}(\beta) \circ_1 \alpha) \stackrel{\mathit{op.comp}}{=}
(\tau_{p+1}(\gamma)\circ_1 \tau_{n+1}(\beta)) \circ_1 \alpha
\stackrel{\eqref{compos_cyclic3}}{=}\\ = \tau_{n+p}(\beta
\circ_{n}\gamma) \circ_1 \alpha \stackrel{\mathit{
\eqref{def_cyclic_compos}}}{=} \alpha\circ_{m+1} (\beta
\circ_{n}\gamma).
\end{multline*}
\end{proof}
\end{Def}

\begin{Rem} In order to avoid unnecessary complications, we are
restricting the discussion to the case where the lowest space
$\mathcal O(1)=k$ is one dimensional; see page \pageref{O(1)=k}. This
means that $\widehat{\mathcal O}(\da,\da;\varnothing)$ is also one
dimensional, and thus there is only one lowest inner product. It
would be an interesting generalization to consider an arbitrary
operad.\\ Another restriction that we chose in order to avoid
complications, is to look at operads $\mathcal O$ which are ungraded.
In the graded case, all the claims of this paper remain true, except
for an adjustment of signs.
\end{Rem}

\subsection{Quadraticity and Koszulness of $\widehat{\mathcal O}$}
\label{quadrat-koszul}

In this subsection we want to derive some implications for
$\widehat{\mathcal O}$, when $\mathcal O$ is quadratic and Koszul.

\begin{Def}[Cyclic Quadratic Operad]
First, recall from \cite{GeK} Section 3, that an operad $\mathcal O$
is called {\em cyclic quadratic}, if it is quadratic, with generators
$E$ and relations $R$, so that the $S_2$-action on $E$ is naturally
extended to a $S_3$-action via the sign-representation, and
$R\subset\mathcal F(E)(3)$ is an $S_4$-invariant subspace. In this
case, $\mathcal O$ becomes a cyclic operad, see \cite{GeK}.
\end{Def}

\begin{Lem}\label{O_hat_quadratic}
Let $\mathcal O$ be cyclic quadratic with generators $E=E^{\ff,\ff}_{\ff}$
and relations $R\subset\mathcal F(E)(3)$. Then $\widehat{\mathcal O}$
is generated by $\widehat{E}$ given by the collection of spaces
\begin{eqnarray*}
\widehat{E} ^{\ff,\ff}_{\ff}:=E^{\ff,\ff}_{\ff}&\subset \widehat {\mathcal
O}(\f,\f;\f),\\%
\widehat{E} ^{\ff,\dda}_{\dda}:=E^{\ff,\ff}_{\ff}&\subset \widehat {\mathcal
O}(\f,\da;\da) ,\\%
\widehat{E} ^{\dda,\ff}_{\dda}:= E^{\ff,\ff}_{\ff}&\subset \widehat {\mathcal
O}(\da,\f;\da),\\%
\widehat{E}^{\dda,\dda}:=k&\subset \widehat {\mathcal O}(\da,\da;\varnothing),
\end{eqnarray*}
and has relations
\begin{eqnarray*}
R\subset\mathcal F(E)(3)\cong\mathcal F(\widehat{E})(\f,\f,\f;\f), \\
R\subset\mathcal F(E)(3)\cong\mathcal F(\widehat{E})(\f,\f,\da;\da), \\
R\subset\mathcal F(E)(3)\cong\mathcal F(\widehat{E})(\f,\da,\f;\da), \\
R\subset\mathcal F(E)(3)\cong\mathcal F(\widehat{E})(\da,\f,\f;\da).
\end{eqnarray*}
together with the relations
\begin{eqnarray*}
G\subset \mathcal F(\widehat{E})(\da,\da,\f;\varnothing), \\ G\subset
\mathcal F(\widehat{E})(\da,\f,\da;\varnothing), \\ G\subset \mathcal
F(\widehat{E})(\f,\da,\da;\varnothing),
\end{eqnarray*}
where $G$ corresponds for a given coloring to the space
\[ G:= span \left<
\begin{pspicture}(0,0.2)(1,1)
 \psline(0.5,0)(0.3,0.5)
 \psline(0.5,0)(0.7,0.5)
 \psline(0.3,0.5)(0.3,1)
 \psline(0.3,0.7)(0.1,1)
 \psline(0.3,0.7)(0.5,1)
 \rput[b](0.7,0.8){$\alpha$}
\end{pspicture}
-
\begin{pspicture}(0,0.2)(2,1)
 \psline(0.5,0)(0.3,0.5)
 \psline(0.5,0)(0.7,0.5)
 \psline(0.7,0.5)(0.7,1)
 \psline(0.7,0.7)(0.9,1)
 \psline(0.7,0.7)(0.5,1)
 \rput[b](1.5,0.6){$\tau_3(\alpha)$}
\end{pspicture}
\text{ , for all } \alpha\in E^{\ff,\ff}_{\ff} \right>.
\]
\end{Lem}
\begin{proof} Lemma~\ref{O_bar_quadratic} guarantees that the
$n$-to-$1$ tree complex modulo relations induces the correct spaces
$\mathcal O(n)$. Thus using the $S_{n+1}$-action on trees with $n+1$
inputs and no output, it remains to show that for $\vec X=(\da,\f,\ldots,
\f,\da)$ one has an identification
\[
\left(\mathcal F(\widehat{E}) /(R,G)\right)(\vec X;\varnothing)
=\widehat{\mathcal O}(\vec X;\varnothing)=\mathcal O(n).
\]
We will do this by looking at the following commutative diagram of
isomorphisms:
\[
\begin{pspicture}(0,0)(10,2.5)
 \rput(2,2){\rnode{A}{$\left(\mathcal F(\widehat{E})/(R,G)\right)
  (\vec X;\varnothing)$}}
 \rput(8,2){\rnode{B}{$\mathcal O(n)$}}
 \rput(8,0){\rnode{C}{$(\mathcal F(E)/(R))(n)$}}
 \ncline{->}{A}{B} \ncput*{$\lambda$}
 \ncline{->}{C}{B} \ncput*{$\mu$}
 \ncline{->}{C}{A} \ncput*{$\kappa$}
\end{pspicture}
\]
First, the map $\kappa$ is induced by the map $\mathcal F(E)\to
\mathcal F(\widehat{E})$, which maps a decorated $n$-to-$1$ binary
tree to the decorated $(n+1)$-to-$0$ binary tree by inverting the
output into a new input by capping it off with an inner product
$1\in\mathcal O(1)=k$, and using the appropriate coloring.
Next, we notice that this map descends to a map $\kappa$, because
$(R)\subset(R,G)$. In fact formula~\eqref{def_cyclic_compos} for
$m=1$ becomes $\alpha\circ_2\beta = \tau_{n+1}(\beta)\circ_1 \alpha$
while $\alpha\in \mathcal O(1)=k$ implies $\tau_{n+1}(\beta)\circ_1
\alpha = \tau_{n+1}(\beta)= \alpha \circ_1 \tau_{n+1}(\beta)$, and
hence
\begin{equation}\label{cyclic_rel}
\alpha\circ_2\beta = \tau_{n+1}(\beta)=\alpha\circ_1
\tau_{n+1}(\beta)
\end{equation}
Notice furthermore that $\kappa $ is surjective, since by the help of
the relations $G$, one can get any element of $\mathcal
F(\widehat{E})$ equivalent to an element in the image of $\kappa$,
i.e., equivalent to a sum of decorated binary trees where the furthest
right input has not branched out at all.\\
Next, we let the map $\mu$ be the isomorphism coming from the
property of $\mathcal O$ being quadratic. To be more precise, $\mu$
evaluates the binary tree according to its decoration using the
composition of the operad.\\
Finally, $\lambda$ is defined similarly to $\mu$, as it evaluates a
binary tree without output according to its decoration, namely
\[
\lambda:(\mathcal F(\widehat{E})/(R,G))(\vec X;\varnothing) \to
\widehat{\mathcal O} (\vec X;\varnothing)=\mathcal O(n).
\]
We claim that $\lambda$ is well defined. After this is shown, the
result follows, since clearly $\lambda\circ \kappa=\mu$ is an
isomorphism, and thus $\kappa$ is injective, and therefore an
isomorphism.\\
Showing that $\lambda$ is well-defined amounts to showing that two
decorated binary trees from $\mathcal F(\widehat{E})$, which are
related by a relation from either $R$ or $G$ induce the same element
in $\mathcal O(n)$ after evaluation. Now, $G$ preserves evaluation,
because equation (\ref{def_cyclic_compos}) for $\alpha\in \mathcal
O(2)= E$ and $\beta=1\in \mathcal O(1)=k$ gives exactly the relations
$G$ and the associativity of the evaluation thus implies the claim.
Similarly, associativity of the composition also implies that the
relations $R$ are preserved, since $\mathcal O$ clearly preserves
those relations.
\end{proof}

\begin{Thm}\label{O_hat_Koszul}
If $\mathcal O$ is cyclic quadratic and Koszul, then
$\widehat{\mathcal O}$ has a resolution given by
\begin{eqnarray*}
\textbf{D}(\widehat{\mathcal O^!}) (\vec X;\f)&\to& \widehat{
\mathcal O}(\vec X;\f)=\mathcal O(n) \\
\textbf{D}(\widehat{\mathcal O^!}) (\vec X;\da)&\to& \widehat{
\mathcal O}(\vec X;\da)=\mathcal O(n) \\
\textbf{D}(\widehat{\mathcal O^!}) (\vec X;\varnothing) &\to&
\widehat{ \mathcal O}(\vec X;\varnothing)=\mathcal O(n-1)
\end{eqnarray*}
where $|\vec X|=n$.
\end{Thm}
The proof of the Theorem will be given in subsection
\ref{proof_main}. For now, we only want to mention that Theorem
\ref{O_hat_Koszul} is the precise version of the following intuitive
combinatorial idea.

\begin{Rem} \label{Koszulness_idea}
Theorem \ref{O_hat_Koszul} states that given a cyclic quadratic
Koszul operad $\mathcal O$, and interpreting outputs as additional
inputs by using the cyclic structure,
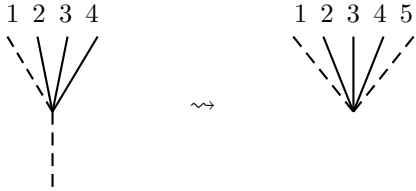
\[
\begin{pspicture}(0,0.6)(4,3.8)
 \psline[linestyle=dashed](2,2)(1.4,3)
 \psline(2,2)(1.8,3)
 \psline(2,2)(2.2,3)
 \psline(2,2)(2.6,3)
 \psline[linestyle=dashed](2,2)(2,1)
 \rput[b](2,3.2){$1\,\,\, 2\,\,\, 3\,\,\, 4$}
 \rput[b](4,2){$\rightsquigarrow$}
\end{pspicture}
\begin{pspicture}(0,0.6)(4,3.8)
 \psline[linestyle=dashed](2,2)(1.2,3)
 \psline(2,2)(1.6,3)
 \psline(2,2)(2,3)
 \psline(2,2)(2.4,3)
 \psline[linestyle=dashed](2,2)(2.8,3)
 \rput[b](2,3.2){$1\,\,\, 2\,\,\, 3\,\,\, 4\,\,\, 5$}
\end{pspicture}
\]
then one can get a resolution of $\mathcal O(n)$ which is now
associated to diagrams with $n+1$ inputs and no output. This
resolution uses a complex of (2-colored) trees with no outputs and is
in fact essentially different to the one provided by Koszulness of
$\mathcal O$, which uses $n$-to-$1$ trees. The reason why Theorem
\ref{O_hat_Koszul} works, strongly uses the fact that
$\widehat{\mathcal O}$ is a 2-colored operad. In fact, up to this
Theorem, there was no need to introduce colors at all. But one can
check that Theorem \ref{O_hat_Koszul} would fail to be true if one
would assume only one color. (In the proof, one uses the two distinct
``dashed'' inputs to think of a diagram as a kind of tensor product
of elements of $\mathcal O^!(k)$, and then tries to mimic the
K\"unneth Theorem.)\\
For example, here is the resolution complex of
$\mathcal O(3)$, interpreted in terms of diagrams without output:
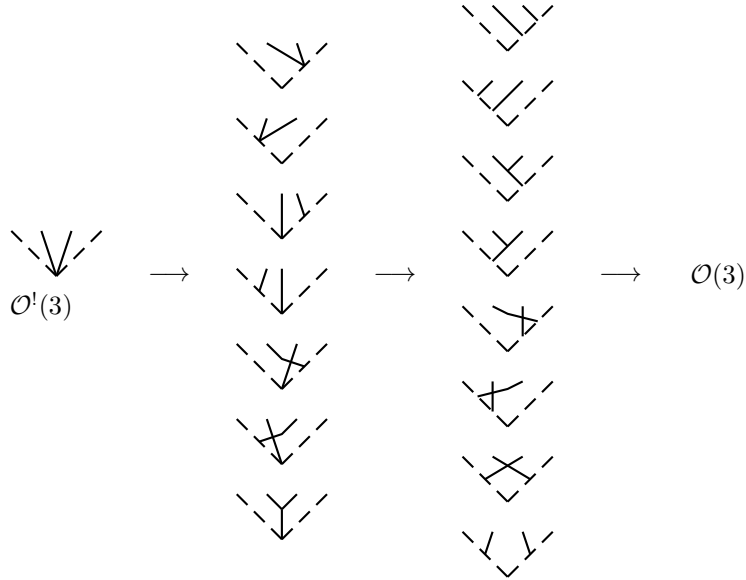
\[
\begin{pspicture}(1,-2)(4,4)
 \psline[linestyle=dashed](2,2)(1.4,2.6)
 \psline[linestyle=dashed](2,2)(2.6,2.6)
 \psline(2,2)(1.8,2.6)
 \psline(2,2)(2.2,2.6)
 \rput(1.8,1.6){$\mathcal O^!(3)$}
 \rput(3.5,2){$\lora$}
\end{pspicture}
\begin{pspicture}(1,-0.5)(3,7)
%
 \psline[linestyle=dashed](2,0)(1.4,0.6)
 \psline[linestyle=dashed](2,0)(2.6,0.6)
 \psline(1.8,0.6)(2,0.4) \psline(2.2,0.6)(2,0.4) \psline(2,0.4)(2,0)
%
 \psline[linestyle=dashed](2,1)(1.4,1.6)
 \psline[linestyle=dashed](2,1)(2.6,1.6)
 \psline(1.8,1.6)(2,1) \psline(2.2,1.6)(2,1.4)(1.7,1.3)
%
 \psline[linestyle=dashed](2,2)(1.4,2.6)
 \psline[linestyle=dashed](2,2)(2.6,2.6)
 \psline(1.8,2.6)(2,2.4)(2.3,2.3) \psline(2.2,2.6)(2,2)
%
 \psline[linestyle=dashed](2,3)(1.4,3.6)
 \psline[linestyle=dashed](2,3)(2.6,3.6)
 \psline(2,3.6)(2,3) \psline(1.8,3.6)(1.7,3.3)
%
 \psline[linestyle=dashed](2,4)(1.4,4.6)
 \psline[linestyle=dashed](2,4)(2.6,4.6)
 \psline(2.2,4.6)(2.3,4.3) \psline(2,4.6)(2,4)
%
 \psline[linestyle=dashed](2,5)(1.4,5.6)
 \psline[linestyle=dashed](2,5)(2.6,5.6)
 \psline(1.8,5.6)(1.7,5.3) \psline(2.2,5.6)(1.7,5.3)
%
 \psline[linestyle=dashed](2,6)(1.4,6.6)
 \psline[linestyle=dashed](2,6)(2.6,6.6)
 \psline(1.8,6.6)(2.3,6.3) \psline(2.2,6.6)(2.3,6.3)
%
\end{pspicture}
\begin{pspicture}(1,-2)(2,4)
 \rput(1.5,2){$\lora$}
\end{pspicture}
\begin{pspicture}(1,0)(3,8)
%
 \psline[linestyle=dashed](2,0)(1.4,0.6)
 \psline[linestyle=dashed](2,0)(2.6,0.6)
 \psline(1.8,0.6)(1.7,0.3) \psline(2.2,0.6)(2.3,0.3)
 \psline[linestyle=dashed](2,1)(1.4,1.6)
 \psline[linestyle=dashed](2,1)(2.6,1.6)
 \psline(1.8,1.6)(2.3,1.3) \psline(2.2,1.6)(1.7,1.3)
 \psline[linestyle=dashed](2,2)(1.4,2.6)
 \psline[linestyle=dashed](2,2)(2.6,2.6)
 \psline(1.8,2.6)(1.8,2.2) \psline(2.2,2.6)(2,2.5)(1.6,2.4)
 \psline[linestyle=dashed](2,3)(1.4,3.6)
 \psline[linestyle=dashed](2,3)(2.6,3.6)
 \psline(1.8,3.6)(2,3.5)(2.4,3.4) \psline(2.2,3.6)(2.2,3.2)
 \psline[linestyle=dashed](2,4)(1.4,4.6)
 \psline[linestyle=dashed](2,4)(2.6,4.6)
 \psline(1.8,4.6)(2,4.4) \psline(2.2,4.6)(1.8,4.2)
 \psline[linestyle=dashed](2,5)(1.4,5.6)
 \psline[linestyle=dashed](2,5)(2.6,5.6)
 \psline(1.8,5.6)(2.2,5.2) \psline(2.2,5.6)(2,5.4)
 \psline[linestyle=dashed](2,6)(1.4,6.6)
 \psline[linestyle=dashed](2,6)(2.6,6.6)
 \psline(1.8,6.6)(1.6,6.4) \psline(2.2,6.6)(1.8,6.2)
\psline[linestyle=dashed](2,7)(1.4,7.6)
 \psline[linestyle=dashed](2,7)(2.6,7.6)
 \psline(1.8,7.6)(2.2,7.2) \psline(2.2,7.6)(2.4,7.4)
\end{pspicture}
\begin{pspicture}(0,-2)(2.5,4)
 \rput(0.5,2){$\lora$}
 \rput(1.8,2){$\mathcal O(3)$}
\end{pspicture}
\]
The dimensions of the individual spaces above are %
$$ dim(\mathcal O^!(3))\to
\begin{array}{c}
  2\cdot dim(\mathcal O^!(3))+\\
  +5\cdot dim(\mathcal O^!(2))^2
\end{array} \to
8\cdot dim(\mathcal O^!(2))^2 \to
 dim(\mathcal O(3)),$$
so that exactness of this complex implies that $$dim(\mathcal
O(3))=3\cdot dim(\mathcal O^!(2))^2-dim(\mathcal O^!(3)),$$ which is
what we expected from the $n$-to-$1$ tree complex of the Koszulness
of $\mathcal O$:
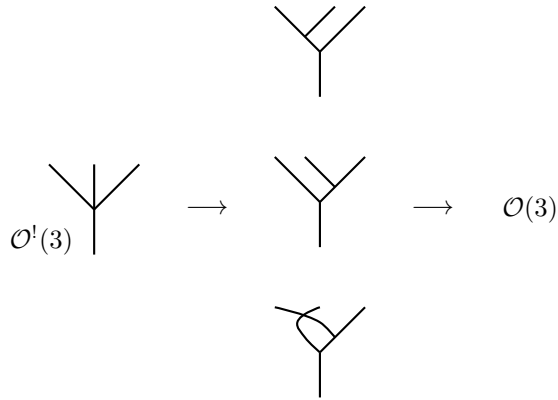
\[
\begin{pspicture}(1,-.5)(4,3.4)
 \psline(2,2)(1.4,2.6)
 \psline(2,2)(2.6,2.6)
 \psline(2,1.4)(2,2.6)
 \rput(1.3,1.6){$\mathcal O^!(3)$}
 \rput(3.5,2){$\lora$}
\end{pspicture}
\begin{pspicture}(1,1.4)(3,6.6)
%
 \pscurve(2,2)(1.8,2.2)(1.7,2.4)(2,2.6) \psline(2,2)(2.6,2.6)
 \pscurve(1.4,2.6)(2,2.4)(2.2,2.2) \psline(2,2)(2,1.4)
 \psline(2,4)(1.4,4.6) \psline(2,4)(2.6,4.6)
 \psline(1.8,4.6)(2.2,4.2) \psline(2,4)(2,3.4)
 \psline(2,6)(1.4,6.6)  \psline(2,6)(2.6,6.6)
 \psline(2.2,6.6)(1.8,6.2) \psline(2,6)(2,5.4)
\end{pspicture}
\begin{pspicture}(0,-.5)(2.5,3.4)
 \rput(0.5,2){$\lora$}
 \rput(1.8,2){$\mathcal O(3)$}
\end{pspicture}
\]
\end{Rem}

\subsection{Inner Products and Homotopy Inner Products}
\label{homotop-ip}

What is the concept of an algebra over $\widehat{\mathcal O}$? Recall
that an algebra over $\widehat{\mathcal O}$ is defined to be a
0/1-operad map $\widehat{\mathcal O}\to \End^{A,M}$ into the
endomorphism 0/1-operad respecting all the structures. Thus, we have
``algebra maps'' $\mathcal O(n) \otimes A^{\otimes n}\to A$ and
``module maps'' $\bigoplus_{r+s=n-1} \mathcal O(n) \otimes A^{\otimes
  r}\otimes M\otimes A^{\otimes s}\to M$. In addition, since
$\widehat{\mathcal O}(\vec X;\varnothing) =\mathcal O(n-1)$ when
$\vec X$ has cardinality $n$ and has precisely two ``dashed'' inputs
(say, in the $i$th and $j$th position), there are also ``inner product
maps''
\[
\mathcal O(n-1) \to Hom(A^{\otimes i-1} \otimes M
\otimes A ^{\otimes j-i-1} \otimes M \otimes A^{\otimes n-j} ,k),
\]
for every $i<j$, depending on the position of the ``dashed'' inputs.
Notice that the lowest case, $n=2$, gives a map $<.,.>:M\otimes M\to k$,
because $\mathcal O(1)=k$ was assumed. Using the composition and
the $S_n$-action of $\widehat{\mathcal O}$, it is easy to see that
all the higher inner product maps are determined by $<.,.>$ together
with the module maps mentioned above. Additionally,
equation~\eqref{cyclic_rel} shows that the inner product satisfies
\[
<\alpha(m_1,a_2,...,a_{n}),m_{n+1}>=<m_1,\tau_{n+1}(\alpha)(a_2,...,
a_{n},m_{n+1})>,
\]
for all $\alpha\in \widehat{\mathcal O}(\da,\f,\ldots,\f;\da)=\mathcal
O(n)$, where we assumed that $m_1, m_{n+1} \in M$ and
$a_2,...,a_{n}\in A$.\\
Our goal is now to find a description for algebras over $\textbf{D}
(\widehat{\mathcal O})$. Clearly, we expect homotopy $\mathcal
O$-algebra and homotopy $\mathcal O$-module maps resolving the algebra
and module maps, respectively, as in
Theorems~\ref{regular_koszuality}~and~\ref{O_bar_algebras}, but we
also expect some sort of ``homotopy inner products'' resolving the
inner product map. For this, we need the definition of a module map
over the modules $F_{\mathcal O,A}M$ and $F_{\mathcal O,A}N$ from
subsection~\ref{ssec:am}.

\begin{Def}[Module Map]\label{modulemap}
A {\em module map} $f\in \mathrm{Mod}(F_{\mathcal O, A}M, F_{\mathcal
O,A}N)$ is a map from $F_{\mathcal O,A}M$ to $F_{\mathcal O,A}N$
making the following diagram commutative:
\[
\begin{CD}
\mathcal O(k) \otimes (F_{\mathcal O}A)^{\otimes (k-1)}\otimes
F_{\mathcal O, A}M & @>\gamma^M>> & F_{\mathcal O, A}M\\
@V\mathrm{id}\otimes \mathrm{id}^{\otimes (k-1)} \otimes f VV & &
@VVfV\\ \mathcal O(k) \otimes (F_{\mathcal O}A)^{\otimes
(k-1)}\otimes F_{\mathcal O, A}N & @>\gamma^N>> & F_{\mathcal O, A}N
\end{CD}
\]
where $\gamma^M$ and $\gamma^N$ are the $(\mathcal O,F_{\mathcal
O}A)$-module structures of $F_{\mathcal O, A} M$ and $F_{\mathcal O,
A} N$ from Definition \ref{freemodule}.\\ Notice that, if $\mathcal
O$ is quadratic, then this condition is equivalent to saying that for
every $\mu\in\mathcal O (2)$, $a\in F_{\mathcal O}A$ and $m\in
F_{\mathcal O, A}M$, one has $ f(\mu(a,m))=(-1)^{|a|\cdot|f|}
\mu(a,f(m))$.\\ The module map $f$ is completely determined by its
restriction to $M\cong \mathcal O(1)\otimes M \subset F_{\mathcal O,
A} M$, i.e., by maps $$f_{n+2}: M \to (\mathcal O(n+1) \otimes
(A^{\otimes n}\otimes N))_{S_{n}\times id},$$ for $n\in\N$.
\end{Def}

We are almost ready to state the Theorem which interprets algebras
over $\textbf{D}(\widehat{\mathcal O})$. The statement of the Theorem
uses the dual structure of a module. In more details, for cyclic
operads $\mathcal O$, we need to use the fact, that if $(A,M)$ is an
algebra over $\textbf{D}(\overline{\mathcal O})$, then so is
$(A,M^*)$, where $M^*$ is the $k$-linear dual of $M$.

\begin{Def}[Dual Homotopy Module] \label{dual-module} Given
differential graded spaced $(A,d_A)$ and $(M,d_M)$, so that $M$ is
finite dimensional in every degree, assume that $M$ is a homotopy
module over a homotopy algebra $A$ over the cyclic operad $\mathcal
O$, i.e., there are given $d\in \mathrm{Der} (F_{\mathcal
O^!}\,sA^*)$ and $g\in \mathrm{Der}_d (F_{\mathcal O^!,sA^*}sM^*)$,
with $d^2=0$ and $g^2=0$. The maps $d$ and $g$ are determined by maps
\begin{eqnarray*}
d_n:\,sA^*&\to& \bigoplus_n \mathcal O^! (n)\otimes (\,sA^*)^{\otimes n},\\ %
g_n:\,sM^*&\to& \bigoplus_{k+l=n-2} \mathcal O^! (k+l+1)\otimes
(\,sA^{*})^{\otimes k} \otimes \,sM^* \otimes (\,sA^{*})^{\otimes l}.
\end{eqnarray*}
Then $M^*$ also becomes a homotopy module over $A$ in the
following way. Define $h\in \mathrm{Der}_d (F_{\mathcal O^!,
sA^*}sM)$ by a map
\begin{equation*}
h_n:\,sM \to \bigoplus_{k+l=n-2} \mathcal O^!(k+l+1)\otimes
(\,sA^{*})^{\otimes k} \otimes\,sM \otimes (\,sA^{*})^{\otimes l},
\end{equation*}
which is given on generators $a_1,...,a_{k},a'_{1},...,a'_l \in \,A
\cong sA$, $m_1\in \,M\cong sM$, $m_2^*\in \,M^*\cong sM^*$ by
\begin{multline*}
h_n(m_1)(a_1,...,a_{k},m_2^*,a'_{1},...,a'_l):=\\ :=(-1)^\epsilon
\cdot \tau_{k+l+2}^{l+1}(g_n(m_2^*)(a'_{1},...,a'_l,m_1,a_1,...,
a_{k})) \in\mathcal O^!(k+l+1).
\end{multline*}
\[
\begin{pspicture}(0,-1)(5,2.5)
 \psline(2,1)(0.8,0)
 \psline(2,1)(1.2,0)
 \psline(2,1)(1.6,0)
 \psline(2,1)(2.4,0)
 \psline(2,1)(2.8,0)
 \psline(2,1)(3.2,0)
 \psline[arrowsize=0.2, arrowinset=0, linestyle=dashed](2,2)(2,1)
 \rput[l](2.5,1.35){$h$} \rput[m](2,2.2){$m_1$}
 \psline[linestyle=dashed](2,0)(2,1)
 \rput(2,-0.4){$a'_{l}$...$a'_{1}$ $m_2^*$ $a_{k}$...$a_{1}$}
 \psline[linearc=0.3, arrowsize=0.15]{->}(1.6,1.2)(1.9,1.5)(2.4,1)(1.9,0.5)(1.5,0.9)
 \rput(4.5,1.5){$:=(-1)^\epsilon \cdot$}
\end{pspicture}
\begin{pspicture}(0,-1)(4,2.5)
 \psline(2,1)(0.8,0)
 \psline(2,1)(1.2,0)
 \psline(2,1)(1.6,0)
 \psline(2,1)(2.4,0)
 \psline(2,1)(2.8,0)
 \psline(2,1)(3.2,0)
 \psline[arrowsize=0.2, arrowinset=0, linestyle=dashed](2,2)(2,1)
 \rput[l](2.5,1.35){$\left(\tau_{k+l+2}^{l+1}\otimes
                      id\right) (g)$} \rput[m](2,2.2){$m_2^*$}
 \psline[linestyle=dashed](2,0)(2,1)
 \rput(2,-0.4){$a_{k}$...$a_{1}$ $m_1$ $a'_{l}$...$a'_{1}$}
 \psline[linearc=0.3, arrowsize=0.15]{->}(1.6,1.2)(1.9,1.5)(2.4,1)(1.9,0.5)(1.5,0.9)
\end{pspicture}
\]
Here $\tau_{k+l+2}^{l+1}$ denotes the $(l+1)$st iteration of
$\tau_{k+l+2}$, and $\epsilon=(|m_1|+|a_1|+...+|a_k|+k+1)\cdot
(|m_2^*|+|a'_1|+...+|a'_l|+l+1)$.
\begin{Lem}
With the above notations, $g^2=0$ implies $h^2=0$.
\end{Lem}
\begin{proof} For better readability, we will show this statement
for the dual maps, namely that $(g^*)^2=0$ implies $(h^*)^2=0$, where
we use the linear dual maps
\begin{eqnarray*}
&d_n^*:&\bigoplus_n \left(\mathcal O^! (n)\right)^* \otimes
(\,s^{-1}A)^{\otimes n}\to \,s^{-1}A,\\ %
&g_n^*:&\bigoplus_{k+l=n-2} \left(\mathcal O^! (k+l+1)\right)^* \otimes
(\,s^{-1}A)^{\otimes k}\otimes \,s^{-1}M \otimes (\,s^{-1}A)^{\otimes l}\to \,s^{-1}M,\\%
&h_n^*:&\bigoplus_{k+l=n-2} \left(\mathcal O^! (k+l+1)\right)^* \otimes
(\,s^{-1}A)^{\otimes k} \otimes \,s^{-1}M^* \otimes (\,s^{-1}A)^{\otimes l}\to \,s^{-1}M^*.
\end{eqnarray*}
Let $a_1,...,a_k,a'_1,...,a'_l\in \,A\cong s^{-1}A$, $m_1^*\in
\,M^*\cong s^{-1}M^*$, $m_2\in \,M\cong s^{-1}M$ and $\alpha^*\in
\mathcal O^!(k+l+1)^*$. Then
$\left((h^*)_n^2(\alpha^*;a_1,...,a_k,m_1^*,a'_1,...,a'_l) \right)
(m_2)$ is given by
\begin{eqnarray*}
&&\sum_{r<s, \beta^*, \gamma^*} (-1)^\zeta
h_p^*(\gamma^*; ..., d_q^*(\beta^*;a_r,...,a_s) ,...,m_1^*,...)(m_2)+\\%
&&\sum_{r, s, \beta^*, \gamma^*} (-1)^\eta h_p^*(\gamma^*;
...,h_q^*(\beta^*;a_r,...,m_1^*,...,a'_s) ,...)(m_2)+\\%
&&\sum_{r<s, \beta^*, \gamma^*} (-1)^\theta
h_p^*(\gamma^*;...,m_1^*,...,d_q^*(\beta^*;a'_r,...,a'_s) ,...)(m_2),
\end{eqnarray*}
where $p$ and $q$ are the appropriate numbers, and the sum in
$\beta^*$ and $\gamma^*$ is the sum in the dual of $\circ_i: \mathcal
O(n)\otimes \mathcal O(m)\to \mathcal O(n+m-1)$, namely
$(\circ_i)^*:\mathcal O(n)^*\to \bigoplus_{r+s=n+1}\mathcal O(r)^*
\otimes \mathcal O(s)^*$, $\alpha^*\mapsto \sum_{i} \gamma^*_i
\otimes \beta^*_i$. The signs $(-1)^\zeta$, $(-1)^\eta$ and
$(-1)^\theta$ are given by the usual rule, that puts $(-1)^{|x|\cdot
|y|}$ whenever a symbol $x$ jumps over a symbol $y$. Now, by
definition, $$h_n^*(\eta^*;a_1,...,a_k, m_1^*,a'_1,...,a'_l)(m_2)=
(-1)^\epsilon g_n^*(\eta^*\circ \tau^{l+1}_{k+l+2};
a'_1,...,a'_l,m_2,a_1,...,a_k)(m_2^*),$$ so that the above expression
is equal to
\begin{eqnarray*}
&&\sum_{r<s, \beta^*, \gamma^*} (-1)^{\zeta'} g_p^*(\gamma^*\circ
\tau_{k+l+2-r+s}^{l+1}; ...,m_2..., d_q^*(\beta^*;a_r,...,a_s) ,...)(m_1^*)+\\%
&&\sum_{r, s, \beta^*, \gamma^*} (-1)^{\eta'} g_p^*(\beta^*\circ
\tau_{k-r+s+3}^{s+1};...,g_q^*(\gamma^*\circ\tau_{l+r-s+1}^{l-s+1};
a_r,...,m_2,...,a'_s) ,...)(m_1^*)+\\%
&&\sum_{r<s, \beta^*, \gamma^*} (-1)^{\theta'} g_p^*(\gamma^*\circ
\tau_{r+s+2-r+s}^{l+1-r+s} ;...,d_q^*(\beta^*;a'_r,...,a'_s)
,...,m_2,...)(m_1^*).
\end{eqnarray*}
But this sum turns out to be nothing but
$(g^*)_n^2(\alpha^*\circ \tau_{k+l+2}^{l+1} ; a'_{1},...,a'_l,
m_2,a_1,...a_{k})(m_1^*)=0,$ because equations \eqref{compos_cyclic2}
and \eqref{compos_cyclic3} show that for any $\beta\in \mathcal
O^!(m)$ and $\gamma\in \mathcal O^!(n)$, it is
\begin{eqnarray*}
\text{for } i\leq n-j:& \tau_{n+m}^{i}(\gamma\circ_j\beta)
= & \tau_{n+1}^{i}(\gamma)\circ_{i+j}\beta, \\%
\text{for } n-j<i\leq n-j+m:& \tau_{n+m}^{i}(\gamma\circ_j\beta)
= & \tau_{m+1}^{i+j-n}(\beta)\circ_{i+j-n}\tau_{n+1}^{n-j}(\gamma), \\%
\text{for } n-j+m<i:& \tau_{n+m}^{i}(\gamma\circ_j\beta) = &
\tau_{n+1}^{i-m+1}(\gamma)\circ_{i+j-m-n}\beta.
\end{eqnarray*}
\end{proof}
\end{Def}

\begin{Thm}\label{O_hat_algebras} Let $(A,d_A)$ and $(M,d_M)$ be
differential graded spaces. Then giving $(A,M)$ the structure of an
algebra over $\textbf{D} (\widehat{\mathcal O})$ is equivalent to the
following data:
\begin{itemize}
\item
a derivation $d\in \mathrm{Der}(F _{\mathcal O}\,sA^*)$ of degree
$1$, with $d^2=0$
\item
a derivation $g\in \mathrm{Der}_d (F_{\mathcal O,\,sA^*}sM^*)$ over
$d$ of degree $1$, with $g^2=0$
\item
a module map $f\in \mathrm{Mod}(F_{\mathcal O, sA^*}sM, F_{\mathcal
O,sA^*}sM^*)$ of degree $0$ such that
\begin{equation}
\label{eq:fhgf}
f\circ h = g \circ f
\end{equation}
and satisfying the following symmetry condition: if $f$ is given by
maps $f_n:\,sM\to \bigoplus_{k+l=n-2} \mathcal O(k+l+1)\otimes
(\,sA^{*})^{\otimes k}\otimes \,sM^* \otimes (\,sA^{*})^{\otimes l}$, then
\begin{multline}
\label{eq:symm}
 f_n(m_2)(\alpha^*;a_1,...,a_i,m_1,a'_1,...,a'_j)=\\
= (-1)^\epsilon f_n(m_1)(\alpha^*\circ \tau_{i+j+2}^{j+1}
;a'_1,...,a'_j,m_2,a_1,...,a_i).
\end{multline}
\end{itemize}
Here, $\epsilon=(|m_2|+|a_1|+...+|a_i|+i+1)\cdot (|m_1|+ |a'_1|+...+
|a'_j|+j+1)$, and $h\in \mathrm{Der}_d (F_{\mathcal O, sA^*}sM)$ is
the homotopy module structure of $M$ induced by $g$, as described in
Definition~\ref{dual-module} above.
\end{Thm}

If $\mathcal O$ is cyclic quadratic and Koszul, then, by Theorem
\ref{O_hat_Koszul}, $\textbf{D}( \widehat{\mathcal O^!} )\cong
\widehat{\mathcal O}$ and we define a {\em homotopy inner product
over} $A$ {\em and} $M$, or {\em homotopy $\mathcal O$-inner
product}, to be given by a module map $f\in \mathrm{Mod}(F_{\mathcal
O^!, sA^*}sM, F_{\mathcal O^!,sA^*}sM^*)$, so that $f\circ h=g \circ
f$.

\begin{proof}
By definition the pair $(A,M)$ has the structure of an algebra over
the 0/1-operad $\mathbf{D}(\widehat {\mathcal O})$ if there exists a
collection of maps of differential graded spaces $$\lambda_{
(\vec X;x)}: \mathbf{D}(\widehat{\mathcal O})(\vec X;x) \to
\End^{A,M} (\vec X;x)$$ for every choice of $\vec X$ and
$x\in\{\f,\da,\varnothing\}$, respecting the composition and the
symmetric group action.\\ In order to simplify our formulas, we set
\begin{align*}
\Delta_n:=&(\f,\ldots,\f),\\
\Gamma_n:=&(\f,\ldots,\f,\da),\\
\Phi_n:=&(\da,\f,\ldots,\f,\da),
\end{align*}
where $n$ is the cardinality of the sequence.\\
The first step in the proof is to establish the following identifications:
\begin{itemize}
\item the $\lambda_{(\Delta_n;\ff)}$'s are equivalent to a derivation
$d\in \mathrm{Der}(F_{\mathcal O}\,sA^*)$;
\item the $\lambda_{(\Gamma_n;\dda)}$'s are equivalent to a derivation
$g\in \mathrm{Der}_d(F_{\mathcal O,\,sA^*}\,sM^*)$;
\item the $\lambda_{(\Phi_n;\varnothing)}$'s are
equivalent to a module map $f\in \mathrm{Mod}
(F_{\mathcal O,\,sA^*}\,sM, F_{\mathcal O,\,sA^*}\,sM^*)$.
\end{itemize}
This is enough since we know from Definition~\ref{def_O_hat} that
these are the only sequences $\vec X$, up to symmetric group action,
which are nontrivial. The cases of having ``dashed'' inputs in the
middle of the sequence $\vec X$ are determined by $\Delta_n$,
$\Gamma_n$, $\Phi_n$ and the symmetric group action.\\ Recall that
$\mathbf{D}(\widehat{\mathcal O})(\vec X;x)$ is defined as
\begin{multline} \label{D(O!)}
   (\widehat {\mathcal O}(\vec X;x))^\vee\stackrel{\de}{\lora}
   \bigoplus_{(\vec X;x)-\mathrm{trees},\,|T|=1}
   (\widehat {\mathcal O}(T))^*\otimes \Det(T) \stackrel{\de}{\lora}
   \dots \\ \dots
   \stackrel{\de}{\lora}
   \bigoplus_{\mathrm{binary}\,\,(\vec X;x)-\mathrm{trees}}
   (\widehat {\mathcal O}(T))^*\otimes \Det(T).
\end{multline}
Now, note that $\lambda$ respecting the composition implies that we
only need to specify $\lambda$ on trees with no internal edge. Then,
any other tree is determined by those trees through composition, and
therefore we also know the action of $\lambda$. Thus, denote by
$\delta_n$, $\gamma_n$ and $\phi_n$ the restrictions of the maps
$\lambda_{( \Delta_n;\ff)}$, $\lambda_{(\Gamma_n;\dda)}$ and
$\lambda_{(\Phi_n;\varnothing)}$ to the left-most components of
$\mathbf{D}(\widehat{\mathcal O})(\Delta_n;\f)$, $\mathbf{D}(
\widehat{\mathcal O})(\Gamma_n;\da)$ and $\mathbf{D}(\widehat{\mathcal
O})(\Phi_n;\varnothing)$ from \eqref{D(O!)} respectively. Explicitly, we
have maps\label{deltagammaphi}
\begin{align*}
\delta_n\colon & \left(\widehat{\mathcal O}(\Delta_n;\f)\right)^\vee
\to (A^*)^{\otimes n}\otimes A,\\ \gamma_n\colon & \left(
\widehat{\mathcal O}(\Gamma_n;\da)\right)^\vee \to (A^*)^{\otimes
(n-1)}\otimes M^*\otimes M,\\ \phi_n\colon & \left(\widehat{\mathcal
O}(\Phi_n;\varnothing)\right)^\vee \to M^*\otimes (A^*)^{\otimes
(n-2)}\otimes M^*.
\end{align*}
Here the $\delta_n$'s are $S_n$-equivariant, the
$\gamma_n$'s are equivariant with respect to the subgroup
$(S_{n-1}\times id)$ of $S_n$ fixing the last component, and the
$\phi_n$'s are equivariant with respect to the subgroup
$(S_{n-2}\times S_2)$ of $S_n$ whose action is given by permuting
the inner $n-2$ components and the outer $2$ components. It is
important to remark that in order to keep track of the twisted sign
in the cobar complex, all these groups act on the left hand side via
the sign representation. Alternatively, we can bring the sign action
to the right hand side and note that acting via the sign
representation on $(A^*)^{\otimes n}$ is like acting via the standard
representation on $(sA^*)^{\otimes n}$. (Similar statements hold for
the other two cases.) Therefore, the above maps become
\begin{align*}
d_n\colon & sA^*\to \left(\mathcal O(n) \otimes (sA^*)^{\otimes n}
\right)_{S_n},\\ g_n\colon & sM^*\to \left(\mathcal O(n)\otimes
(sA^*)^{\otimes (n-1)}\otimes sM^*\right)_{S_{n-1}\times id},\\
f_n\colon & sM\to \left(\mathcal O(n-1) \otimes (sA^*)^{\otimes
(n-2)} \otimes sM^*\right)_{S_{n-2}\times id},
\end{align*}
where the $f_n$ have the additional property of respecting the
$S_2$-action of switching the $M$-factors. This is precisely
described by the symmetry condition \eqref{eq:symm}.

We now want to use the fact that $\lambda$ is a map of differential
graded operads, and thus the maps $\lambda_{(\vec X;x)}$ have to
respect the differentials. We know that the differential $\partial$
in the cobar dual decollapses the edges, whereas on $\End^{A,M}$ we
have the differential $\widetilde{\partial}^{A,M}$ induced by the
differentials of the differential graded vector spaces $(A,d_A)$ and
$(M,d_M)$.\\ Observe first that $\lambda_{(\vec X;\ff)}$ respecting the
differentials means
\[
\begin{CD}
\left(\widehat{\mathcal O}(\Delta_n;\f)\right)^\vee & @>\partial>> &
\bigoplus_{|T|=1} \left(\widehat{\mathcal O}(\Delta_{k_T};\f)
\right)^\vee \otimes \left(\widehat{\mathcal O} (\Delta_{l_T};\f)
\right)^\vee\\ @V\delta_n VV & & @VV ev_T \circ( \delta_{k_T}\otimes
\delta_{l_T})V\\ Hom(A^{\otimes n},A)& @>\widetilde{
\partial}^{A,M}>> & Hom(A^{\otimes n},A),
\end{CD}
\]
where we look at trees $T$ with exactly one internal edge and we thus
get the induced maps $\delta_{k_T}$ and $\delta_{l_T}$ from splitting
$T$ into two trees with no internal edges and $k_T$, respectively
$l_T$, external edges. The map $ev_T$ is the composition of maps
$(ev_T(\mu \otimes \nu)) (a_1,...,a_r)=\mu(a_1,...,\nu(...),...,a_r)$
according to $T$. Here $\widetilde{\partial}^{A,M}$ is induced by the
differential $d_A$ of the differential graded vector space $(A,d_A)$
as follows:
\begin{align*}
(\widetilde{\partial}^{A,M}(\chi))(a_1,\ldots,a_n) =
d_A(\chi(a_1,\ldots,a_n)) + \sum_{i=1}^n (-1)^{\epsilon_{i}}
\chi(\ldots,d_A(a_i),\ldots).
\end{align*}
where $(-1)^{\epsilon_{i}}$ is the usual sign, see
page~\pageref{sgn_epsilon_i}. Thus, for $n\ge 2$ we have
\begin{equation}\label{delcircdel}
\widetilde{\partial}^{A,M}\circ \delta_n=\sum_{|T|=1}
ev_T\circ(\delta_{k_T}\otimes \delta_{l_T} )\circ\partial.
\end{equation}
Let $\delta_1:\widehat{\mathcal O}(\f;\f)^\vee=\mathcal O(1)^\vee\to
A^*\otimes A$ be defined by $d_A:A\to A$ and denote its associated
dual map by $d_1:sA^*\to \mathcal O(1)\otimes sA^*$. Now extend all
$d_i$, for $i\ge 1$, to derivations $d_i\in \mathrm{Der} (F_{\mathcal
O} sA^*)$, then equation \eqref{delcircdel} becomes $d_1\circ
d_n+d_n\circ d_1=\sum_{j=2}^{n-1} d_j\circ d_{n-j+1}$. Thus, if $d:=
\sum_{i\ge 1} d_i$, then equation \eqref{delcircdel} is equivalent to
$d^2=0$.

In a similar way, $\lambda_{(\vec X;\dda)}$ respecting the differentials
is equivalent to a derivation $g\in \mathrm{Der}_d (F_{\mathcal
O,sA^*}sM^*)$ over $d$ with $g^2=0$. Here, the lowest component $g_1:
sM^* \to sM^*=\mathcal O(1)\otimes sM^*$ is induced by the
differential $d_M$ on $M$.

We claim that $\lambda_{(\Phi_n;\varnothing)}$ respecting the
differentials is equivalent to equation~\eqref{eq:fhgf}. Again, we
start from the commutative diagram
\[
\begin{CD}
\left(\widehat{\mathcal O}(\Phi_n;\varnothing)\right)^\vee & @>\partial>> &
\bigoplus_{|T|=1} \left(\widehat{\mathcal O}(\Lambda_{k_T};l) \right)^\vee
\otimes \left(\widehat{\mathcal O}
(\Phi_{l_T};\varnothing) \right)^\vee\\ @V\phi_n VV & & @VV \oplus_T
(\mathit{ev}_T \circ(\lambda_{k_T}\otimes \phi_{l_T}))V\\
\Hom(M\otimes A^{\otimes n-2}\otimes M,k)&
@>\widetilde{ \partial}^{A,M}>> & \Hom(M\otimes A^{\otimes n-2}\otimes
M,k),
\end{CD}
\]
where $\Lambda_{k_T}$ is either $\Delta_{k_T}$ or $\Gamma_{k_T}$,
$\lambda_{k_T}$ is either $\delta_{k_T}$ or $\gamma_{k_T}$, $l$ is
either ``full'' or ``dashed'' depending on the coloring of the
internal edge of the tree $T$. Also recall that $\widetilde{
\partial}^{A,M}$ is the lift of the differentials $d_A$ and $d_M$ to
an $\End^{A,M}$-derivation. We thus get the equation
\begin{equation}\label{other}
\widetilde{\partial}^{A,M}\circ \phi_n=\sum_{|T|=1} ev_T\circ
(\lambda_{k_T} \otimes \phi_{l_T} )\circ\partial,
\end{equation}
where the left-hand-side is given by
\begin{eqnarray*}
(\widetilde{\partial}^{A,M}\circ \phi_n(\alpha))(m_1,a_1,...,
a_{n-2},m_2) &=& (\phi_n (\alpha))(d_M(m_1),a_1,...,a_{n-2},m_2)+\\
&&+\sum_i (-1)^{\epsilon_{i-1}} (\phi_n
(\alpha))(m_1,...,d_A(a_i),...,m_2)\\ &&(-1)^{\epsilon_{n}}(\phi_n
(\alpha))(m_1,a_1,...,a_{n-2},d_M(m_2)).
\end{eqnarray*}
with, again, the signs defined on page~\pageref{sgn_epsilon_i}. The
right-hand-side of \eqref{other} is all combinations of plugging
trees into an inner product diagram. Thus, if we take the partial sum
of all trees that are composed by one of the first $n-1$ inputs of
the inner product diagram, together with all but the last summand of
the left-hand-side of \eqref{other}, then this readily corresponds to
the term $g\circ f$.\\ We now show that the composition at the last
input of the inner product diagram corresponds to $f\circ h$. In
order to do this, it is easier to use the dual maps
\begin{eqnarray*}
&f^*&:\bigoplus_{n}(\mathcal O(n+1))^* \otimes s^{-1}M \otimes
s^{-1}A^{\otimes n}\to s^{-1}M^*\\%
&g^*&:\bigoplus_{n}(\mathcal O(n+1))^*\otimes s^{-1}A^{\otimes n}
\otimes s^{-1}M\to s^{-1}M\\%
&h^*&:\bigoplus_{n}(\mathcal O(n+1))^*\otimes s^{-1}M^* \otimes
s^{-1}A^{\otimes n}\to s^{-1}M^*,
\end{eqnarray*}
where $h$ is the induced map coming from Definition
\ref{dual-module}. Now, $(h^* \circ f^*)( \alpha^*;m,a_1,...,a_n)\in
s^{-1}M^*$ and thus $((h^*\circ f^*)( \alpha^* ;m,a_1,...,a_n))(m')$
is equal to
\begin{eqnarray*}
&&\sum_{k,\alpha_1^*, \alpha_2^*} (h^*(\alpha_1^*;f^*(\alpha_2^*;m,
a_1,...,a_k),a_{k+1}, ...,a_n ))(m')=\\%
&=&\sum_{k,\alpha_1^*, \alpha_2^*}(-1)^\epsilon
(g^*(\alpha_1^*\circ\tau^{-1}_{n-k+2}; a_{k+1}, ..., a_n,m'))
(f^*(\alpha_2^*;m,a_1,...,a_k))=\\%
&=&\sum_{k,\alpha_1^*, \alpha_2^*}(-1)^{\epsilon'} (f^*(\alpha_2^*;
m,a_1,...,a_k)) (g^*(\alpha_1^*\circ\tau^{-1}_{n-k+2};a_{k+1},
...,a_n,m'))=\\%
&=&\sum_{k,\alpha_1^*, \alpha_2^*} (\phi(\alpha_2^*)) \left(m,a_1,
...,a_k, (\gamma(\alpha_1^* \circ \tau^{-1}_{n-k+2})) (a_{k+1},
...,a_n,m')\right),
\end{eqnarray*}
where the signs come from Definition \ref{dual-module}, and $\phi$
and $\gamma$ are induced by the maps $\phi_n$ and $\gamma_n$ from page
\pageref{deltagammaphi}. The above sum is over all $\alpha_1^*$ and
$\alpha_2^*$ which appear in the dual of the composition
$\alpha_1\circ_1\alpha_2=\alpha$, namely $(\circ_1)^*:\alpha^*\mapsto
\sum_i \alpha_{1,i}^{*} \otimes \alpha_{2,i}^{*}$. But using equation
\eqref{def_cyclic_compos}, it follows that $\alpha_1\circ_1
\alpha_2=\alpha$ is equivalent to $\alpha_2\circ_{k+2} \tau^{-1}_{
n-k+2}(\alpha_1 )=\alpha$. Thus, we get exactly all combinations of
splitting the inner product diagram at the last spot.
\end{proof}

\subsection{Examples: $\mathcal Assoc$, $\mathcal Comm$, $\mathcal Lie$}
\label{Ex:A,C,L}

\begin{Exa}[$\mathcal Assoc$]\label{assoc}
Define the associative operad $\mathcal Assoc$ to be the cyclic
quadratic operad with the $2$-dimensional generators $E=span(S_2)$
with the usual $S_2$-action on $E$ extended to an $S_3$-action by the
sign-action. If $A$ is an algebra over $\mathcal Assoc$, then the
elements of $E$ are thought of as multiplications $\mu(a,b)=a\cdot b$
and $\mu^{op}(a,b)=b\cdot a$ in $A$. In this setting, the defining
relations of $\mathcal Assoc$ correspond to the $6$-dimensional
subspace of relations $(a_{\sigma(1)}\cdot a_{\sigma(2)})\cdot
a_{\sigma(3)}=a_{\sigma(1)}\cdot (a_{\sigma(2)}\cdot a_{\sigma(3)})$,
for all $\sigma\in S_3$, of the $12$-dimensional $\mathcal F(E)(3)$.
It can be shown that this defines a cyclic quadratic Koszul operad,
which is its own quadratic dual $\mathcal Assoc^!=\mathcal Assoc$,
see \cite{GeK}.

More explicitly, one can show that the $n$th space of $\mathcal
Assoc$ is given by $\mathcal Assoc(n)=$ $k$-vector space generated by
$S_n$. In order to distinguish between $\sigma\in S_n$ and a basis of
$\mathcal Assoc(n)$, we use the notation $[\sigma]$ for a basis of
$\mathcal Assoc(n)= span\{[\sigma]\}_{\sigma\in S_n}$. Before
describing the $S_{n+1}$-action and the composition, we want to
mention that the correspondence of algebras $A$ over $\mathcal
Assoc$ and associative algebras is given by letting the algebra
map $\mathcal Assoc(n)= span\{ [\sigma]\}_{\sigma\in S_n} \to
Hom(A^{\otimes n},A)$ be defined as $[\sigma]:(a_1,...,a_n)\mapsto
a_{\sigma(1)}\cdot ...\cdot a_{\sigma(n)}$.
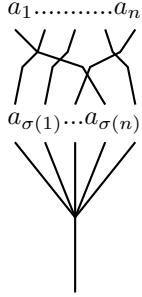
\[
\begin{pspicture}(0,0.3)(4,4.6) 
 \psline(2,1.5)(1.2,2.5)  \psline(1.2,3)(1.3,3.5)(1.4,3.6)(1.5,3.7)(1.6,4)
 \psline(2,1.5)(1.6,2.5)
 \psline(1.6,3)(1.7,3.5)(1.8,3.6)(1.9,3.7)(2,4)
 \psline(2,1.5)(2,2.5)
 \psline(2,3)(2.2,3.5)(2.4,3.6)(2.6,3.7)(2.8,4)
 \psline(2,1.5)(2.4,2.5)  \psline(2.4,3)(2.1,3.5)(1.8,3.6)(1.5,3.7)(1.2,4)
 \psline(2,1.5)(2.8,2.5)  \psline(2.8,3)(2.7,3.5)(2.6,3.6)(2.5,3.7)(2.4,4)
 \psline(2,1.5)(2,0.5)
 \rput(2,4.25){$a_1...........a_n$}
 \rput(2,2.75){$a_{\sigma(1)}...a_{\sigma(n)}$}
\end{pspicture}
\]
The $S_{n+1}$-action on $\mathcal Assoc(n)$ is given in the following
way. As before, denote by $\tau_{n+1}\in S_{n+1}$ the cycle $(12 ...
n+1)$. If $\mu\in S_{n+1}$, $\sigma\in S_n\subset S_{n+1}$, then
notice that one can rewrite $\mu\cdot\sigma$ as a product
$\mu\cdot\sigma =\sigma'\cdot (\tau_{n+1}) ^k$, where both
$\sigma'\in S_n\subset S_{n+1}$ and the power $k\in \{0,...,n\}$ are
uniquely determined. In this notation, we define the action by
$\mu.[\sigma]:=[\sigma']$. This implies for example that
$\tau_{n+1}.[id_{\{1,...,n\}}]= [id_{\{1,...,n\}}]$, and if $\mu\in
S_n$, $\sigma\in S_n$, then $\mu.[\sigma]=[\mu\cdot \sigma]$.\\%
Next, we define the composition $\circ_i:\mathcal Assoc(n) \otimes
\mathcal Assoc(m)\to \mathcal Assoc(n+m-1)$. If $\mu\in S_n$ and
$\sigma\in S_m$, then we set
$$ [\mu]\circ_i[\sigma]:=[(id\times ...\times id\times \sigma\times
id...\times id)\cdot \mu_{1,...,1,m,1,...,1}], $$%
where $\times:S_k\times S_l\to S_{k+l}$ adjoins two permutations and
$\mu_{1,...,1,m,1,...,1}\in S_{n+m-1}$ permutes blocks of sizes $1$,
..., $1$, $m$, $1$, ..., $1$. Thus, $\sigma$ is plugged into $\mu$ into
the $i$-th slot.

We want to describe the homotopy $\mathcal Assoc$-algebras, -modules
and -inner products in more details. The data required for a homotopy
$\mathcal Assoc$-algebra are maps $\mathcal Assoc(n)^\vee \otimes
A^{\otimes n}\to A$, which respect the $S_n$-action in the domain.
Using the explicit description of $\mathcal Assoc(n)=span(S_n)$, we
see that this determines maps $d_n:sA^*\to (sA^*)^{\otimes n}$, for
$n\geq 2$. The map $d:=d_1+d_2+...$, where $d_1:sA^*\to sA^*$ comes
from the differential on $A$, has to satisfy $d^2=0$, when lifted to
the tensor algebra $T(sA^*)$ as a derivation.\\%
Similarly a homotopy $\mathcal Assoc$-module $M$ over $A$ is given by
maps $\bigoplus_{p+q\geq 1} \mathcal Assoc(p+q+1)^\vee \otimes
A^{\otimes p}\otimes M \otimes A^{\otimes q}\to M$, which determine
the maps $$ g_{p,q}:sM^*\to sA^{* \otimes p}\otimes sM^* \otimes
sA^{* \otimes q}, $$ after factoring out the $S_n$-action. The map
$g:=\sum_{p,q} g_{p,q}$, where $g_{0,0}:sM^*\to sM^*$ comes from the
differential on $M$, has to be a differential, i.e., $g^2=0$.\\%
Finally, a homotopy $\mathcal Assoc$-inner product consists of maps
$\bigoplus_{p+q\geq 1} \mathcal Assoc(p+q+1)^\vee \otimes A^{\otimes
p}\otimes M \otimes A^{\otimes q}\to M^*,$ which reduce to maps $$
f_{p,q}\in sA^{*\otimes p}\otimes sM^* \otimes sA^{* \otimes q}
\otimes sM^*. $$ The map $f:=\sum_{r,s}f_{r,s}$ satisfies $f\circ
g=h\circ f$, where $h$ is the dual homotopy $\mathcal Assoc$-module
map induced by $g$, see
Definition \ref{dual-module}.\\%
Thus, we see that the above concept determines an $\infty$-inner
product over an $A_\infty$-algebra as defined in \cite{Tr}, which
additionally satisfies the symmetry condition \eqref{eq:symm}, coming
from switching the $M$ factors. Furthermore, one can see from
\cite{Tr2} Lemma 2.14, that this additional symmetry \eqref{eq:symm}
implies that the $\infty$-inner product is symmetric in the sense of
\cite{Tr2} Definition 2.13.
\end{Exa}

\begin{Rem}
From the above example, we see that $\widehat{\mathcal
Assoc}$-algebras determine symmetric $\infty$-inner products. In
\cite{Tr2}, it was shown that non-degenerate, symmetric
$\infty$-inner products induce a BV-structure on the
Hochschild-cohomology of the given A$_\infty$ algebra. It would be
interesting to have a generalization of this result to any cyclic
operad $\mathcal O$, which amounts to look for a similar ``$\mathcal
O$-BV-structure'' on the homology of the chain complex $
C_n^{\mathcal O}(A)=\left(\mathcal O^!(n)\otimes A^{\otimes n}
\right)_{S_n} $ of an $\mathcal O$-algebra $A$ from \cite{GK},
(4.2.1)-(4.2.4).
\end{Rem}

\begin{Exa}[$\mathcal Comm$, $\mathcal Lie$]
Let $\mathcal Comm$ be the quadratic operad defined by the
1-dimensional generators $E=k$ and relations $R\subset\mathcal
F(E)(3)$ identifying all trees with each other. It follows that
$\mathcal Comm(n)=k$ is a one-dimensional space for each $n$. The
$S_{n+1}$-action, $\sigma:k\lora k$, and the compositions, $k\circ_i
k\lora k$, are the trivial maps, and it is clear that this defines a
cyclic quadratic operad. Define $\mathcal Lie$ to be the quadratic
dual of $\mathcal Comm$, i.e., $\mathcal Lie=\mathcal Comm^!$ and
$\mathcal Lie^!= \mathcal Comm$. It can be shown that $\mathcal
Comm$, and thus also $\mathcal Lie$, are Koszul (see \cite{GK}).\\
What are the homotopy concepts for $\mathcal Comm$ and $\mathcal
Lie$? A homotopy $\mathcal Lie$-algebra consists of maps $\mathcal
Comm(n)^\vee\otimes A^{\otimes n}\to A$, which respect the symmetric
group action in the domain. As $\mathcal Comm(n)=k$, this amounts to
specifying a derivation into the symmetric algebra $S(sA^*)$
generated by $sA^*$, i.e., $d:sA^*\to S(sA^*)$, with $d^2=0$. A
homotopy $\mathcal Lie$-inner product over $A$, where $A$ is
interpreted as a module over itself, requires maps $\mathcal
Comm(n+1)^\vee\otimes A^{\otimes n}\otimes A\to A^*$, or more
explicitly, maps $<...>_n:A^{\otimes n+2}\to k$ for each $n\geq 0$.
Thus homotopy $\mathcal Lie$-inner products have exactly one higher
inner product map $<...>_n$ in each dimension, satisfying certain
relations.\\ Similarly, one can show that homotopy $\mathcal
Comm$-algebras are specified by a derivation $d: sA^*\to L(sA^*)$,
where $L(sA^*)$ is the free Lie-algebra generated by $sA^*$, with
$d^2=0$. Homotopy $\mathcal Comm$-inner products require maps
$<...>_x:A^{\otimes n+2}\to k$, where the set of labels $x$
corresponds to a set of generators of $\mathcal Lie(n+1)$.
\end{Exa}

\begin{Rem} A more explicit approach to homotopy inner products
over $\mathcal Lie$ and $\mathcal Comm$ was developed by Scott Wilson
\cite{W} independently of the approach given here. It seems that the
two concepts of homotopy inner products do in fact coincide.
\end{Rem}

As an application of the above structure, we want to show how
homotopy $\mathcal Comm$-inner products can be obtained on the chain
level of a topological space. The construction for the homotopy
$\mathcal Comm$ algebra is taken from D. Sullivan's paper \cite{S} on
the construction of local infinity structures, and the one for the
homotopy $\mathcal Comm$ inner product is a slight variation of the
same idea. We will only sketch the main ideas of the proof since the
full details were already given for the associative case by M.
Zeinalian and the second author in \cite{TZ}.

\begin{Prop}\label{prop:comm-pd}
Let $C$ be a finite dimensional simplicial complex so
that the closure of every simplex is contractible, and suppose we are
given a fundamental cycle $\mu$ in $C$. Then $A:=C^*$ has the structure
of a $\widehat{\mathcal Comm}$ algebra, such that the lowest
multiplication is the symmetrized Alexander-Whitney multiplication
and the lowest inner product is given by capping with $\mu$.
\end{Prop}
\begin{proof}
First, let's construct the homotopy $\mathcal Comm$ algebra on $C^*$.
Let $L(sC)=L_1\oplus L_2\oplus...$ be the free Lie algebra on $sC$,
decomposed by monomial degree in $sC$. One needs to find a derivation
$d:L(sC)\to L(sC)$ of degree $1$, given by maps $d=d_1+d_2+...$,
where $d_i:sC\to L_i$, so that $d^2=0$. Let $d_1:sC\to sC$ be the
differential on $C$, and $d_2$ be the symmetrized Alexander-Whitney
comultiplication. For the general $d_i$, we use the inductive
hypothesis that $d_1$, ..., $d_{i-1}$ are local maps so that
$\nabla_i:=d_1+ ...+d_{i-1}$ has a square $\nabla_i^2$ mapping only
into higher components $L_{i}$, $L_{i+1}$,.... Here, ``local'' means
that every simplex maps into the sub-Lie algebra of its closure. Now,
by the Jacobi-identity, it is $0=[\nabla_i, [\nabla_i,\nabla_i]]=
[d_1,e_i]\text{+ higher terms}$, where $e_i:sC\to L_i$ is the lowest
term of $[\nabla_i,\nabla_i]$. Thus $e_i$ is $[d_1,.]$-closed and
thus, using the hypothesis of the proposition, also locally
$[d_1,.]$-exact. These local terms can be put together to give a map
$d_i$, so that $[d_1,d_i]$ vanishes on  $L_1\oplus ... \oplus
L_{i-1}$ and equals $- 1/2\cdot e_i$ on $L_i$. In other words, $(d_1+
...+d_i)^2=1/2\cdot [d_1+ ...+d_i,d_1+ ...+d_i]=1/2\cdot [\nabla_i
,\nabla_i ]+ [d_1,d_i]+\text{higher terms}$, maps only into
$L_{i+1}\oplus L_{i+2}\oplus...$. This completes the inductive step,
and thus produces the wanted homotopy $\mathcal Comm$ algebra. But
then any homotopy $\mathcal Comm$ algebra also induces a homotopy
$\mathcal Comm$ module over itself.\\ Finally, we want to construct
the homotopy $\mathcal Comm$ inner product. For this, first recall
that a homotopy $\mathcal Comm$ inner product consists of an element
$f$ in $\mathcal M:= \mathrm{Mod}(F_{\mathcal Comm, sC}sC^*,
F_{\mathcal Comm ,sC}sC)$ of degree $0$, so that $f\circ h=g\circ f$,
where $g$ and $h$ are the induced homotopy module structures for
$C^*$ and $C$ coming from the homotopy algebra structure $d$. The
space $\mathcal M$ can now be endowed with the differential
$\mathcal D(f):=f\circ h-(-1)^{|f|} g\circ f$, satisfying
$\mathcal D^2=0$, since both $g^2=0$ and $h^2=0$. By Definition
\ref{modulemap}, we can identify $\mathcal M$ with the space
\begin{equation*} \mathcal M\cong \bigoplus_{i\geq 0} \Hom\left(
sC^*,(\mathcal Comm (i+1)\otimes (sC^{\otimes i}\otimes sC)
)_{S_i\times id}\right)\cong\bigoplus_{i\geq 0} \mathcal M_i,
\end{equation*}
where $\mathcal M_i:=\left(C\otimes C^{\otimes i}\otimes C
\right)_{id\times S_i\times id}$ is the component of monomial degree
$i+2$, and where we chose to ignore the internal degree for
simplicity. Notice that the differential $\mathcal D$ maps $\mathcal
M_i$ to an equal or higher monomial degree $\bigoplus_{j\geq
i}\mathcal M_j$, with lowest component being given by applying $d_1$.
Furthermore $\mathcal D$ is a local map, since $d$ is a local map.\\
Our goal is then to construct a chain map $\chi: (C,d_1) \to
(\mathcal M,\mathcal D)$ by induction on the monomial degree of
$\chi=\chi_0+\chi_1+\chi_2+...$, where $\chi_i:C\to \mathcal M_i$.
Define $\chi_0:C\to C\otimes C$ to be the symmetrized
Alexander-Whitney comultiplication. Now assume inductively, that one
has local maps $\chi_0$, ..., $\chi_{i-1}$, so that $\Upsilon_i:=
\chi_0+...+\chi_{i-1}$ commutes with $d_1$ and $\mathcal D$ up to
monomial degree greater or equal to $i$: $\Upsilon_i \circ d_1-
\mathcal D\circ \Upsilon_i :C\to \mathcal M_i\oplus \mathcal M_{i+1}
\oplus...$. Then $0=\left(. \circ d_1- \mathcal D\circ .\right) \circ
\left(. \circ d_1- \mathcal D\circ .\right) (\Upsilon_i)=\left(.
\circ d_1- d_1\circ .\right)\left( \Upsilon_i \circ d_1- \mathcal
D\circ \Upsilon_i\right)+\text{higher terms}= \epsilon_i \circ
d_1-d_1\circ \epsilon_i+\text{higher terms}$, where $\epsilon_i:C \to
\mathcal M_i$ is the lowest term of $\Upsilon_i \circ d_1- \mathcal
D\circ \Upsilon_i$. Thus $\epsilon_i$ is $(.\circ d_1-d_1\circ
.)$-closed and thus by the hypothesis of the proposition also locally
$(.\circ d_1-d_1\circ .)$-exact. Those local terms can be put
together to give a map $\chi_i:C\to \mathcal M_i$ with $\chi_i\circ
d_1-d_1 \circ \chi_i=-\epsilon_i+\text{higher terms}$. In other
words, one has $(\chi_0+...+\chi_i) \circ d_1- \mathcal D\circ
(\chi_0+...+ \chi_i)= (\chi_0+...+\chi_{i-1}) \circ d_1- \mathcal
D\circ (\chi_0+...+ \chi_{i-1})+(\chi_i\circ d_1-\mathcal D\circ
\chi_i)=\epsilon_i-\epsilon_i+\text{higher terms}=\text{higher
terms}:C\to \mathcal M_{i+1} \oplus \mathcal M_{i+2} \oplus...$. This
completes the inductive step.\\
Having the chain map $\chi:C\to \mathcal M$, we use the fundamental
cycle $\mu\in C$ in order to define $f:=\chi(\mu)$. Clearly
$f\circ h-g\circ f=\mathcal
D(f)=\mathcal D(\chi(\mu))=\chi(d_1( \mu))=0$, which shows that $f$
satisfies equation \eqref{eq:fhgf}. Furthermore, the symmetry
condition \eqref{eq:symm} can be shown analogously to \cite{TZ},
Proposition 3.6.
\end{proof}

\subsection{Proof of Theorem \ref{O_hat_Koszul}}
\label{proof_main}

\setcounter{Thm}{7}
\begin{Thm}
If $\mathcal O$ is cyclic quadratic and Koszul, then
$\widehat{\mathcal O}$ has a resolution given by
\begin{eqnarray*}
\textbf{D}(\widehat{\mathcal O^!}) (\vec X;\f)&\to& \widehat{
\mathcal O}(\vec X;\f)=\mathcal O(n) \\%
\textbf{D}(\widehat{\mathcal O^!}) (\vec X;\da)&\to& \widehat{
\mathcal O}(\vec X;\da)=\mathcal O(n) \\%
\textbf{D}(\widehat{\mathcal O^!}) (\vec X;\varnothing)&\to& \widehat{
\mathcal O}(\vec X;\varnothing)=\mathcal O(n-1)
\end{eqnarray*}
where $|\vec X|=n$.
\end{Thm}

\begin{proof} The first and second quasi-isomorphisms follow from
Proposition~\ref{O_bar_Koszul}. As for the lower one, we need to show
that the homology of $\textbf{D}(\widehat{\mathcal O^!}) (\vec X;
\varnothing)$ is concentrated in degree $0$:
\begin{eqnarray}\label{H0}
H_0 \left(\textbf{D}(\widehat{\mathcal O^!}) (\vec X;\varnothing)\right
)&=&\widehat{ \mathcal O}(\vec X;\varnothing)\\ %
\label{Hi<0} H_{i<0} \left(\textbf{D}(\widehat{\mathcal O^!})
(\vec X;\varnothing)\right)&=&\{0\}
\end{eqnarray}
It is again enough to restrict attention to the case where $\vec X =
(\da,\f,\ldots,\f,\da)$. Recall that the map $\textbf{D} (\widehat{\mathcal
O^!}) (\vec X;\varnothing)\to \widehat{ \mathcal
O}(\vec X;\varnothing)=\mathcal O(n-1)$ ends with the spaces
(cf. equation~\eqref{sequence_H0})
\begin{equation*}
\cdots \stackrel{\de}{\lora} \textbf{D}(\widehat{\mathcal O^!})
(\vec X;\varnothing)^{-1} \stackrel{\de}{\lora}
\textbf{D}(\widehat{\mathcal O^!}) (\vec X;\varnothing)^0
\stackrel{proj}\lora \widehat{\mathcal O}(\vec X;\varnothing).
\end{equation*}
As $\mathcal O$ and thus $\mathcal O^!$ are quadratic, we have the
following identification, using the language and results of Lemma
\ref{O_hat_quadratic}:
\begin{eqnarray*}
\widehat{\mathcal O}(\vec X;\varnothing)&=&\mathcal F(
\widehat{E})/(R,G) (\vec X;\varnothing) \\
\textbf{D}(\widehat{\mathcal O^!}) (\vec X;\varnothing)^0&=&
\bigoplus_{ \substack{
  \text{binary trees }T \\
  \text{of type }(\vec X;\varnothing)
}} (\widehat{E^\vee}(T))\otimes \Det(T) =\mathcal
F(\widehat{E})(\vec X;\varnothing)\\ \textbf{D}(\widehat{\mathcal
O^!}) (\vec X;\varnothing)^{-1}&=&\bigoplus_{ \substack{
  \text{trees }T\text{ of type }(A(f,d;\varnothing)), \\
  \text{binary vertices, except}\\
  \text{one ternary vertex}
}} \left(\widehat{\mathcal O^!}(T)\right)\otimes \Det(T) =\\ &=&
\left\{
\begin{array}{c}
\text{space of relations in $\mathcal F(\widehat{E})(\vec X;\varnothing)$}\\
\text{generated by $R$ and $G$}
\end{array}
\right\}
\end{eqnarray*}
The last equality follows, because the inner product relations for
$\mathcal O^!$ (namely, the relation space $G$ for the cyclic
quadratic operad $\mathcal O^!$ from Lemma \ref{O_hat_quadratic})
are the orthogonal complement of the inner product relations for
$\mathcal O$. Hence it is clear that the map $proj$ is surjective with
kernel $\de\left(\textbf{D}(\widehat{\mathcal O^!}) (\vec X;
\varnothing )^{-1}\right)$. We thus showed equation (\ref{H0}).

As for equation (\ref{Hi<0}), we will use an induction that shows
that every closed element in $\textbf{D}(\widehat{\mathcal O^!})
(\vec X;\varnothing)^{-r}$, for $r\geq 1$ is also exact. Intuitively,
fixing the two ``dashed'' inputs can be used to think of $\textbf{D}
(\widehat{\mathcal O^!}) (\vec X;\varnothing)^{-r}$ as a (twisted)
tensor product of exact complexes. We then want to apply a version of
the K\"unneth Theorem to show that $\textbf{D} (\widehat{\mathcal
O^!}) (\vec X;\varnothing)^{-r}$ is exact. Practically, this amounts
to sliding all of the ``full'' inputs from one of the two ``dashed''
inputs to the other.\\ We need the following definition. Given two
decorated trees $\varphi\in \textbf{D}(\mathcal O^!)(k), \psi \in
\textbf{D}(\mathcal O^!)(l)$, $k+l=n$, we define the operations $*$
and $\#$, so that both $\varphi * \psi$ and $\varphi \# \psi$ belong to
$\textbf{D}(\widehat{ \mathcal O^!})(\vec X;\varnothing)$, in the
following way. First, for $\varphi * \psi$ take the outputs of
$\varphi$ and $\psi$ and insert them into the unique inner product
decorated by $1\in \widehat{\mathcal O^!}(\da,\da;\varnothing)$:
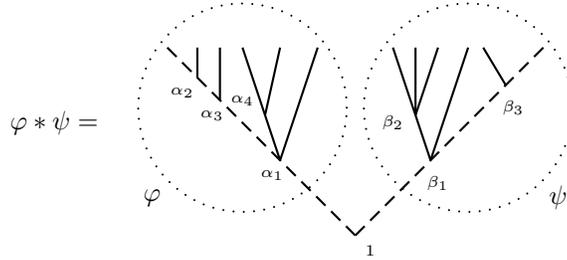
\[
\begin{pspicture}(-3,0)(6,4)
 \psline[linestyle=dashed](2.5,0.5)(0,3)
 \psline[linestyle=dashed](2.5,0.5)(5,3)
 \rput(2.7,0.3){\tiny $1$}
 \pscircle[linestyle=dotted](4,2.2){1.4}
 \rput(-0.2,1){$\varphi$}
 \psline(1.5,1.5)(2,3)
 \psline(1.5,1.5)(1,3)
 \psline(1.3,2.1)(1.5,3)
 \psline(0.4,2.6)(0.4,3)
 \psline(0.7,2.3)(0.7,3)
 \rput(1.4,1.3){\tiny $\alpha_1$}
 \rput(0.2,2.4){\tiny $\alpha_2$}
 \rput(0.6,2.1){\tiny $\alpha_3$}
 \rput(1,2.3){\tiny $\alpha_4$}
 \pscircle[linestyle=dotted](1,2.2){1.4}
 \rput(5.2,1){$\psi$}
 \psline(3.5,1.5)(3,3)
 \psline(3.5,1.5)(4,3)
 \psline(4.5,2.5)(4.2,3)
 \psline(3.3,2.1)(3.3,3)
 \psline(3.3,2.1)(3.6,3)
 \rput(3.6,1.2){\tiny $\beta_1$}
 \rput(3,2){\tiny $\beta_2$}
 \rput(4.6,2.2){\tiny $\beta_3$}
 \rput(-1.5,2){$\varphi * \psi =$}
\end{pspicture}
\]
The operation $\#$ is defined slightly differently. Suppose that
$\varphi\in \textbf{D}(\mathcal O^!)(k)$, where $k\geq 2$. Then, one
first identifies $\varphi$ with an element in $\textbf{D} (\widehat{
\mathcal O^!})(\da,\f,\ldots,\f,\da;\varnothing)$ by interpreting the
lowest decoration $\alpha_1\in O^!(m)$ of $\varphi$, as an inner
product $\alpha_1\in \widehat{\mathcal
O^!}(\da,\f,\ldots,\f,\da;\varnothing)=\mathcal O^!(m)$ rather than a
multiplication $\alpha_1\in\widehat{\mathcal
O^!}(\da,\f,\ldots,\f;\da)=\mathcal O^!(m)$. Then attach $\psi$ to the
input on the furthest right:
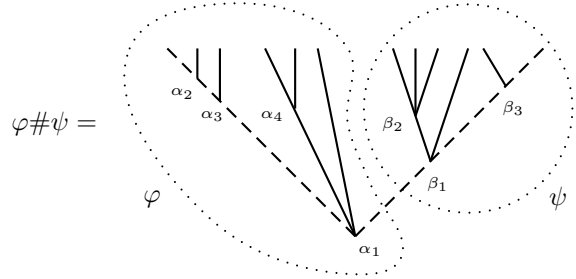
\[
\begin{pspicture}(-3,0)(6,4)
 \psline[linestyle=dashed](2.5,0.5)(0,3)
 \psline[linestyle=dashed](2.5,0.5)(5,3)
 \rput(2.7,0.3){\tiny $\alpha_1$}
 \pscircle[linestyle=dotted](4,2.2){1.4}
 \rput(-0.2,1){$\varphi$}
 \psline(2.5,0.5)(2,3)
 \psline(2.5,0.5)(1.3,3)
 \psline(1.7,2.2)(1.7,3)
 \psline(0.4,2.6)(0.4,3)
 \psline(0.7,2.3)(0.7,3)
 \rput(0.2,2.4){\tiny $\alpha_2$}
 \rput(0.6,2.1){\tiny $\alpha_3$}
 \rput(1.4,2.1){\tiny $\alpha_4$}
 \psccurve[linestyle=dotted](0.5,0.8)(-0.5,3)(1,3.6)(2.4,3)(2.4,2)(2.7,1.2)(3.1,0.2)(2.4,0)
 \rput(5.2,1){$\psi$}
 \psline(3.5,1.5)(3,3)
 \psline(3.5,1.5)(4,3)
 \psline(4.5,2.5)(4.2,3)
 \psline(3.3,2.1)(3.3,3)
 \psline(3.3,2.1)(3.6,3)
 \rput(3.6,1.2){\tiny $\beta_1$}
 \rput(3,2){\tiny $\beta_2$}
 \rput(4.6,2.2){\tiny $\beta_3$}
 \rput(-1.5,2){$\varphi \# \psi =$}
\end{pspicture}
\]
Both $\varphi$ and $\psi$ are elements of $\textbf{D} (\mathcal O^!)$
and thus uncolored. The coloration for $*$ and $\#$ is uniquely
determined by having the first and last entry dashed. After defining
$*$ and $\#$ on generators of $\textbf{D} (\mathcal O^!)$, we extend
them bilinearily to maps $\textbf{D} (\mathcal O^!)\otimes \textbf{D}
(\mathcal O^!)\to \textbf{D} (\widehat{\mathcal O^!})(\vec X; \varnothing)$.

Recall that a $(k,l)$-shuffling of two ordered set of indices $i_1<
\ldots <i_k$ and $j_1< \ldots <j_l$ is a permutation $\sigma$ of
$\{i_1, \ldots i_k, j_1, \ldots,j_l\}$ such that $\sigma(i_1)<\ldots
\sigma(i_k)$ and $\sigma(j_1)<\ldots \sigma(j_l)$. Suppose that we
have a labeled tree $\varphi$ with $k+1$ inputs and a labeled tree
$\psi$ with $l+1$ inputs. As mentioned before, we want to restrict
our attention to the case of elements in $\textbf{D}
(\widehat{\mathcal O^!})(\vec X; \varnothing)$ whose ``dashed''
inputs are labeled to be the first and the last input, and thus
appear in the planar representation on the far left and the far
right. (The other cases are analogous or alternatively can be deduced
from this one.) Then for every $(k,l)$-shuffling $\sigma$, define
$\phi *_\sigma \psi$, resp. $\phi \#_\sigma \psi$, as the composition
of $*$, resp. $\#$, with $\sigma$ applied to the ``full'' leaves of
the resulting labeled tree. The ``dashed'' inputs remain far left and
far right. The importance of this definition lies in the fact that
every labeled tree in $\textbf{D}(\widehat {\mathcal
O^!})(\vec X;\varnothing)$, whose first and last inputs are
``dashed'', can uniquely be written in the form $\varphi
*_\sigma\psi$ or $\varphi\#_\sigma\psi$.

The argument is now to perform an induction with the following inductive
statement for $s\in \mathbb N$:
\begin{itemize}
\item[]
Let $\chi\in \textbf{D}(\widehat{\mathcal O^!})
(\vec X;\varnothing)^{-r}$ be a closed element $\de(\chi)=0$. Then
$\chi$ is homologous to a sum $\sum_i \sum_\sigma \varphi_i *_\sigma
\psi_i+\sum_j \sum_\sigma \varphi'_j \#_\sigma\psi'_j$, where the
total degree of each $\psi_i$ and each $\psi'_j$ is less or equal to
$-s$, i.e., they are elements of $\textbf{D}(\mathcal O^!)(l=l(\psi)
)^{degree\leq(-s)}$.
\end{itemize}
Intuitively, the induction will change $\chi$ homologously to
decorated trees whose total degree is more and more concentrated on
the right branch of the tree.\\
The case $s=1$ will be shown in Lemma \ref{s=1}, whereas the
induction step can be found in Lemma \ref{s->s+1}. The proof of the
Theorem is then completed by noticing that $\textbf{D}(\mathcal O^!)
(l)$ is concentrated in finite degrees, so that the $\psi_i$'s and
$\psi'_j$'s eventually have to be $0$. In other words, every closed
element $\chi$ is homologous to $0$, which means that it is exact.
Thus the complex $ \textbf{D}(\widehat{\mathcal O^!})
(\vec X;\varnothing)$ has no homology in degrees $r\neq 0$.
\end{proof}

\setcounter{Thm}{18}

We are interested in the compatibility of $*$ and $\#$ with the
differentials. In order to make this explicit, we need to brake the
differential into several parts.
\begin{equation*}
\cdots \lora \textbf{D}(\widehat{\mathcal O^!}) (\vec X;\varnothing)^{-r-1}
\stackrel{\de}{\lora} \textbf{D}(\widehat{\mathcal O^!})
(\vec X;\varnothing)^{-r}\lora \cdots
\end{equation*}
The action of $\de$ is given by taking a decorated tree and to expand
one of its (non-binary) vertices. We write the differential%
$$ \de=\de_{l}+\de_{m}+\de_{r}+
\de_{il}+\de_{im}+\de_{ir} $$%
according to the location of this expanded vertex and the location of
its expansion.\\ If the expanded vertex is not the lowest ``inner
product vertex'', then use the first three differentials. The
expanded vertex lies either on a branch which is branched off the
left ``dashed'' input (use $\de_l$), or it lies on a branch which
leads into the inner product (use $\de_m$), or it lies on a branch
which is branched off the right ``dashed'' input (use $\de_r$):
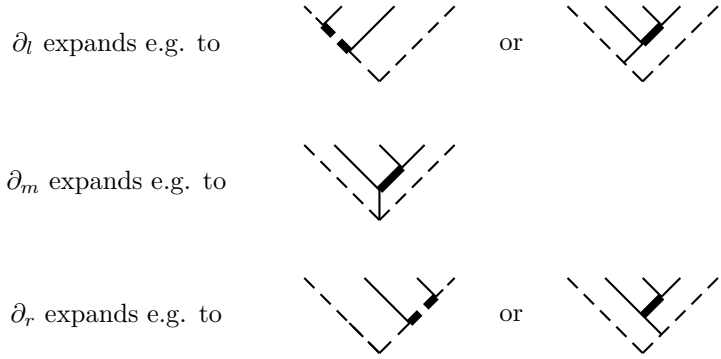
\[
\begin{pspicture}(-3,1.2)(8.5,2.7)
 \psline[linestyle=dashed](1,2.5)(1.25,2.25)
 \psline[linewidth=3pt, linestyle=dashed](1.25,2.25)(1.6,1.9)
 \psline[linestyle=dashed](1.6,1.9)(2,1.5)
 \psline(1.25,2.25)(1.5,2.5)
 \psline(1.6,1.9)(2.2,2.5)
 \psline[linestyle=dashed](3,2.5)(2,1.5)
 \rput(-1.5,2){$\de_l$ expands e.g. to}
 \rput(3.75,2){or}
 \psline[linestyle=dashed](4.5,2.5)(5.5,1.5)
 \psline[linestyle=dashed](6.5,2.5)(5.5,1.5)
 \psline(5.25,1.75)(5.5,2)
 \psline(5.5,2)(5,2.5)
 \psline[linewidth=3pt](5.5,2)(5.75,2.25)
 \psline(5.75,2.25)(6,2.5)
 \psline(5.75,2.25)(5.5,2.5)
\end{pspicture}
\]
\[
\begin{pspicture}(-3,1.2)(8.5,2.7)
 \psline[linestyle=dashed](1,2.5)(2,1.5)
 \psline(2,1.5)(2,1.9)
 \psline(2,1.9)(2.6,2.5)
 \psline(2,1.9)(1.4,2.5)
 \psline(2.3,2.2)(2,2.5)
 \psline[linewidth=3pt](2,1.9)(2.3,2.2)
 \psline[linestyle=dashed](3,2.5)(2,1.5)
 \rput(-1.5,2){$\de_m$ expands e.g. to}
\end{pspicture}
\]
\[
\begin{pspicture}(-3,1.2)(8.5,2.7)
 \psline[linestyle=dashed](1,2.5)(2,1.5)
 \psline[linewidth=3pt, linestyle=dashed](2.75,2.25)(2.4,1.9)
 \psline[linestyle=dashed](1.6,1.9)(2,1.5)
 \psline(2.75,2.25)(2.5,2.5)
 \psline(2.4,1.9)(1.8,2.5)
 \psline[linestyle=dashed](3,2.5)(2.75,2.25)
 \psline[linestyle=dashed](2.4,1.9)(2,1.5)
 \rput(-1.5,2){$\de_r$ expands e.g. to}
 \rput(3.75,2){or}
 \psline[linestyle=dashed](4.5,2.5)(5.5,1.5)
 \psline[linestyle=dashed](6.5,2.5)(5.5,1.5)
 \psline(5.75,1.75)(5,2.5)
 \psline[linewidth=3pt](5.5,2)(5.75,2.25)
 \psline(5.75,2.25)(6,2.5) \psline(5.75,2.25)(5.5,2.5)
\end{pspicture}
\]
Here, the fat lines indicate the expanded vertices.\\ %
If the ``inner product vertex'' is expanded, then use one of the last
three differentials:
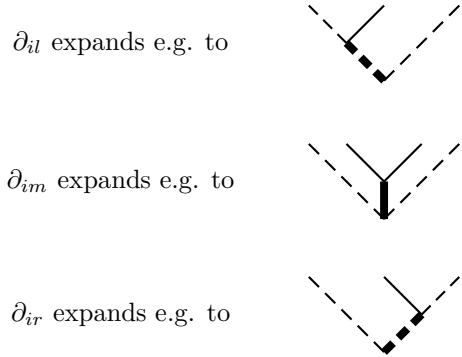
\[
\begin{pspicture}(-3,1.2)(4.5,2.7)
 \psline[linestyle=dashed](1,2.5)(1.5,2)
 \psline[linewidth=3pt, linestyle=dashed](1.5,2)(2,1.5)
 \psline(1.5,2)(2,2.5)
 \psline[linestyle=dashed](3,2.5)(2,1.5)
 \rput(-1.5,2){$\de_{il}$ expands e.g. to}
\end{pspicture}
\]
\[
\begin{pspicture}(-3,1.2)(4.5,2.7)
 \psline[linestyle=dashed](1,2.5)(2,1.5)
 \psline[linewidth=3pt](2,1.5)(2,2)
 \psline(2,2)(1.5,2.5)
 \psline(2,2)(2.5,2.5)
 \psline[linestyle=dashed](3,2.5)(2,1.5)
 \rput(-1.5,2){$\de_{im}$ expands e.g. to}
\end{pspicture}
\]
\[
\begin{pspicture}(-3,1.2)(4.5,2.7)
 \psline[linestyle=dashed](1,2.5)(2,1.5)
 \psline[linewidth=3pt, linestyle=dashed](2,1.5)(2.5,2)
 \psline[linestyle=dashed](2.5,2)(3,2.5)
 \psline(2.5,2)(2,2.5)
 \rput(-1.5,2){$\de_{ir}$ expands e.g. to}
\end{pspicture}
\]
Notice that
\begin{equation}\label{delt(*)}
\de(\varphi * \psi)= \de_l(\varphi * \psi)+ \de_r(\varphi
* \psi)=\de(\varphi) * \psi +\varphi * \de (\psi),
\end{equation}
so that $*$ is a chain map, whereas the compatibility of $\#$
with $\de$ is given by
\begin{equation}\label{delt(**)}
\de(\varphi \# \psi)=\de(\varphi) \# \psi +\varphi \#
\de(\psi)\pm\varphi * \psi\pm \de_{ir}(\varphi \# \psi).
\end{equation}
In order to obtain the correct signs, we assumed the identification
$\Det(T_1 \circ_{k+2} T_2)\cong \Det(T_1)\otimes \Det(T_2)$ in the
definition of $*$ and $\#$.\\ Clearly, the same equations
hold for $*_\sigma$ and $\#_\sigma$ since the shuffling is applied on
the top of the trees and does not interact with $\de$.

\begin{Lem}\label{s=1}
Let $\chi\in \textbf{D}(\widehat{\mathcal O^!})
(\vec X;\varnothing)^{-r}$ be a closed element $\de(\chi)=0$. Then
$\chi$ is homologous to a sum $\sum_i \sum _\sigma \varphi_i *_\sigma
\psi_i+\sum_j \sum _\sigma \varphi'_j \#_\sigma\psi'_j$, where no
$\psi_i$ or $\psi'_j$ is of total degree $0$, i.e.,  no $\psi_i$ or
$\psi'_j$ is a binary tree.
\end{Lem}
\begin{proof} We prove this Lemma recursively on the number of
inputs of $\psi_i$ and $\psi'_j$, namely, we prove that the following
statement by induction on $t$:
\begin{itemize}
\item[]
Let $\chi\in \textbf{D}(\widehat{\mathcal O^!})
(\vec X;\varnothing)^{-r}$ be a closed element $\de(\chi)=0$. Then
$\chi$ is homologous to a sum $\sum_i \sum_\sigma \varphi_i *_\sigma
\psi_i + \sum_j \sum_\sigma \varphi'_j \#_\sigma\psi'_j$, where none
of the $\psi_i$'s or the $\psi'_j$'s are binary with less or equal to
$t$ inputs. In other words, none of the $\psi$'s are elements of
$\textbf{D}(\mathcal O^!) (l) ^0$ for $l\leq t$.
\end{itemize}
We get the claim of the Lemma simply by taking $t=n$, so that we
eliminated all $\psi$'s of degree $0$.

$[t=1]$: We need to eliminate terms of the form $\tilde{\varphi}_i
*_\sigma 1$ and $\tilde{\varphi}'_j \#_\sigma 1$ from $\chi=\sum_i
\sum_\sigma \varphi_i *_\sigma \psi_i+\sum_j \sum_\sigma \varphi'_j
\#_\sigma \psi'_j$, where $1\in \textbf{D}(\mathcal O^!) (1)$. Let's
start with $\tilde{\varphi}'_j \#_\sigma 1$. As the $\de(\tilde{
\varphi}'_j) \#_\sigma 1$ are the only terms in $\de(\chi)=0$ which
are of the form $\varphi \#_\sigma 1$ (see equation
\eqref{delt(**)}), it follows that $\de(\sum_j \tilde{
\varphi}'_j)=0$. Thus one can find an element $\Phi'\in
\textbf{D}(\mathcal O^!)$ with $\de(\Phi') = \sum_j\tilde{
\varphi}'_j$. Then subtract $\de(\Phi'\#_\sigma 1)= \de(\Phi')
\#_\sigma 1+\Phi' *_\sigma 1+\de_{ir}(\Phi' \#_\sigma 1)$ from
$\chi$. The result, which we call again $\chi$, won't contain any
terms of the form $\varphi\#_\sigma 1$ but it possibly introduced
more ($\varphi *_\sigma 1$)-terms. Next, we can take the new $\chi$
and do a similar argument to eliminate the ($\varphi *_\sigma
1$)-terms. Namely, using the fact that there are no ($\varphi
\#_\sigma 1$)-terms, it follows from (\ref{delt(*)}) that the
$\de(\tilde{\varphi}_i) *_\sigma 1$ are the only terms in $\de(\chi)
=0$ of the form $\varphi *_\sigma 1$. Thus it must be
$\de(\sum_i\tilde{\varphi}_i)=0$, so that there exists an element
$\Phi\in \textbf{D}(\mathcal O^!)$ with $\de(\Phi) = \sum_i
\tilde{\varphi}_i$. After subtracting $\de(\Phi *_\sigma 1)=\de(\Phi)
*_\sigma 1$ from $\chi$, all the ($\varphi *_\sigma 1$)-terms
disappeared without having introduced any ($\varphi \#_\sigma
1$)-terms.

$[t-1\to t]$: We use the same idea in the induction step that we used
for the induction start. Namely, let $\chi=\sum_i \sum_\sigma
\varphi_i *_\sigma \psi_i+\sum_j \sum_\sigma \varphi'_j \#_\sigma
\psi'_j$, where none of the $\psi$'s is both binary and has strictly
less than $t$ inputs. Assume that in the above sum for $\chi$, the
terms $\sum_i \sum_\sigma \tilde{\varphi}_i *_\sigma \tilde{\psi}_i +
\sum_j \sum_\sigma \tilde{\varphi}'_j \#_\sigma \tilde{\psi}'_j$ are
the ones where the $\psi$'s have exactly $t$ inputs, i.e.,
$\tilde{\psi}_i$, and $\tilde{\psi}'_j \in\textbf{D}(\mathcal O^!)
(t)$. Our goal is to eliminate those $\tilde{\psi}$'s of degree $0$.
Again, we first want to focus on the sum $\sum_j \sum_\sigma
\tilde{\varphi}'_j \#_\sigma \tilde{\psi}'_j$. We claim that the
tensor product $\sum_j \tilde{\varphi}'_j \otimes \tilde{\psi}'_j \in
\textbf{D} (\mathcal O^!) (n-t)\otimes \textbf{D}(\mathcal O^!)
(t)$ has the property%
$$ \de\left( \sum_j\tilde{ \varphi}'_j \otimes \tilde{\psi}'_j
\right)  \in \bigoplus_{l\leq -1}\textbf{D}(\mathcal O^!)
(n-t)^k\otimes \textbf{D} (\mathcal O^!) (t)^l , $$%
namely, its differential in the complex of the tensor product has no
$\psi$-component in degree $0$. Again the reason for this comes from
the fact that $\de(\chi)=0$, and the fact that the only terms in
$\de(\chi)$, which have binary $\psi$'s with $t$ inputs are of the
form $\de(\tilde{\varphi}'_j)\#_\sigma\tilde{ \psi}'_j$ or
$\tilde{\varphi}'_j \#_\sigma \de(\tilde{\psi}'_j)$ or are terms in
the image of $*_\sigma$. Thus we can use Lemma~\ref{Kuenneth} below
to see that there exists an element $\sum_q \tilde{\Phi}'_q \otimes
\tilde{\Psi}'_q$ with%
$$ \sum_j \tilde{\varphi}'_j \otimes \tilde{\psi}'_j -\de\left(
\sum_q \tilde{\Phi}'_q \otimes \tilde{\Psi}'_q\right)\in
\bigoplus_{l\leq -1} \textbf{D}(\mathcal O^!) (n-t)^{k}\otimes
\textbf{D}(\mathcal O^!) (t)^l . $$%
If we take $\chi^0 :=\chi-\de(\sum_q \sum_\sigma
\tilde{\Phi}'_q \#_\sigma \tilde{\Psi}'_q)\in
\textbf{D}(\widehat{\mathcal O^!}) (\vec X;\varnothing)$, then the
terms $\tilde{\varphi}'_j\#_\sigma \tilde{\psi}'_j$, whose $\tilde{
\psi}$'s have $t$ inputs, cancel. Therefore, the only binary $\psi$'s with
$t$ inputs are in the image of the map $*_\sigma$. Write them as
$\sum_i \sum_\sigma \tilde{\tilde{\varphi}}_i *_\sigma \tilde{
\tilde{\psi}}_i$. Then $\de(\chi^0)=0$ implies that $\sum_i
\tilde{\tilde{ \varphi}}_i \otimes \tilde{\tilde{\psi}}_i$ is closed
and thus exact in the tensor complex $\textbf{D}(\mathcal O^!)
(n-t)\otimes \textbf{D} (\mathcal O^!) (t)$. If we write $\de(\sum_q
\tilde{\tilde{ \Phi}}_q \otimes \tilde{\tilde {\Psi}}_q)=\sum_i
\tilde{\tilde{ \varphi}}_i \otimes \tilde{ \tilde{\psi}}_i$, then
$\chi^0-\de( \sum_q \sum_\sigma \tilde{\tilde{ \Phi}}_q *_\sigma
\tilde{\tilde{\Psi}}_q)$ has no more binary $\psi$'s with $t$ inputs.
\end{proof}

The proof of the following elementary Lemma, which was used in
Lemma \ref{s=1}, is left to the reader.
\begin{Lem}\label{Kuenneth} Let $C_*$ and $D_*$ be two complexes
concentrated in non-positive degrees, with homology concentrated in
degree $0$. Let $\rho\in C_*\otimes D_*$ be an element of non-zero
total degree. Then, if $\rho$ is closed up to degree $r\leq -1$ in
$D_*$, i.e.,
$$ \de(\rho)\in \bigoplus_{l\leq r} C_k\otimes D_l, $$
then it is also exact up to degree $r$ in $D_r$, i.e.,
$$ \exists\ \phi\in C_*\otimes D_*:\quad  \rho-\de(\phi) \in
\bigoplus_{l\leq r} C_k\otimes D_l. $$
\end{Lem}

\begin{Lem}\label{s->s+1}
Let $\chi=\sum_i \sum_\sigma \varphi_i *_\sigma \psi_i+\sum_j
\sum_\sigma \varphi'_j \#_\sigma\psi'_j \in\textbf{D}
(\widehat{\mathcal O^!}) (\vec X;\varnothing)^{-r}$ be a closed
element $\de(\chi)=0$, so that the total degree of each $\psi_i$ and
each $\psi'_j$ is less or equal to $-s$. Then $\chi$ is homologous to
a similar sum, where the total degree of the $\psi$'s is less or
equal to $-s-1$.
\end{Lem}
\begin{proof}
This Lemma uses the same idea as the above Lemma \ref{s=1}, but turns
out to be easier to perform. Namely, let $\chi=\sum_i \sum_\sigma
\varphi_i *_\sigma \psi_i+\sum_j\sum_\sigma  \varphi'_j
\#_\sigma\psi'_j$, where we assume an expansion so that
$\{\varphi_i\}_i$ are linear independent, $\{\varphi'_j \}_j$ are
linear independent, but the $\psi_i$ and the $\psi'_j$ are allowed to
be linear combinations in $\textbf{D}(\mathcal O^!)$. We claim that
those elements $\psi_i$ and $\psi'_j$, which are of degree $-s$, are
closed in $\textbf{D} (\mathcal O^!)$. This follows from
$\de(\chi)=0$ and the inductive hypothesis, because the only terms of
$\de(\chi)$, whose $\psi$'s are of degree $-s+1$, are the terms
$\de_r(\varphi_i *_\sigma \psi_i)=\varphi_i *_\sigma \de(\psi_i)$ and
$\de_r(\varphi'_j \#_\sigma \psi'_j)=\varphi'_j \#_\sigma
\de(\psi'_j)$, where $\psi_i$ and $\psi'_j$ are necessarily of degree
$-s$. The exactness of $\textbf{D} (\mathcal O^!)$ at the degree $-s$
implies that the degree $-s$ $\psi$'s are exact:
$\psi_i=\de(\Psi_i)$, $\psi'_j =\de(\Psi'_j)$. Now, take
$S:=\de\left( \sum_i\sum_\sigma \varphi_i *_\sigma \Psi_i+\sum_j
\sum_\sigma \varphi'_j \#_\sigma\Psi'_j\right)$, where the sum is
over those $i$'s and $j$'s which have the constructed $\Psi_i$'s and
$\Psi'_j$'s. As the total degree of the $\Psi$'s is $-s-1$, we see
that the only terms of $S$ with degree greater or equal to $-s$ are
the terms
\begin{multline*}
 \de_r\left( \sum_i \sum_\sigma \varphi_i *_\sigma \Psi_i+\sum_j
\sum_\sigma \varphi'_j \#_\sigma\Psi'_j\right)=\\
=\sum_i \sum_\sigma \varphi_i *_\sigma \de(\Psi_i) +
\sum_j\sum_\sigma \varphi'_j \#_\sigma\de(\Psi'_j)=\\
=\sum_i \sum_\sigma \varphi_i *_\sigma \psi_i+
\sum_j \sum_\sigma \varphi'_j \#_\sigma\psi'_j
\end{multline*}
It follows that $\chi-S$ only contains terms of degree less or equal
to $-s-1$, and $\chi$ is homologous to $\chi-S$.
\end{proof}



\begin{thebibliography}{99}

\bibitem[A]{Ad} J. F. Adams, ``Infinite loop space'' Princeton
University Press (1978)

\bibitem[BV]{BV} J. M. Boardman, R. M. Vogt, ``Homotopy Invariant
Algebraic Structures on Topological Spaces'' \lnm{347}, Springer
Verlag (1973)

\bibitem[G]{Gan} W. L. Gan, ``Koszul duality for dioperads'',
\mrl{10}, no. 1, 109-124 (2003)

\bibitem[GJ1]{GJ} E. Getzler, J. D. S. Jones, ``$A_\infty$-algebras and
the cyclic bar complex'' \ijm{34}, no. 2, 256-283 (1990)

\bibitem[GJ2]{GJ2} E. Getzler, J. D. S. Jones, ``Operads, homotopy
algebra and the iterated integrals for the double loop spaces''
\texttt{hep-th/9403055}

\bibitem[GeK]{GeK} E. Getzler, M. M. Kapranov, ``Cyclic Operads and Cyclic
Homology'' Geometry, Topology, and Physics for Raoul Bott,
International Press (1994)

\bibitem[GiK]{GK} V. Ginzburg, M. M. Kapranov, ``Koszul duality for
operads'' \dmj{76}, no. 1, 203-272 (1994)

\bibitem[L1]{L1} P. van der Laan, ``Operads, Hopf algebras and coloured
Koszul duality'' Utrecht University Ph.D. Thesis (2003)

\bibitem[L2]{L2} P. van der Laan, ``Coloured Koszul duality and strongly
homotopy operads'' \texttt{math.QA/0312147}

\bibitem[M1]{M} M. Markl, ``Homotopy Algebras are Homotopy Algebras''
\texttt{math.AT/9907138}

\bibitem[M2]{M2} M. Markl, ``Homotopy Diagrams of Algebras''
Proceedings of the 21st Winter School ``Geometry and Physics'',
Rend.\ Circ.\ Mat.\ Palermo (2) Suppl.\ no. 69, 161--180 (2002)

\bibitem[MSS]{MSS} M. Markl, S. Shnider, J. Stasheff,  ``Operads in
Algebra, Topology and Physics'' AMS Mathematical Surveys and
Monographs {\bf 96} (2002)

\bibitem[S]{S} D. Sullivan, ``Local constructions of infinity structures'' to
appear

\bibitem[T1]{Tr} T. Tradler, ``Infinity-inner-products on
A-infinity-algebras'' \texttt{math.AT/010827}

\bibitem[T2]{Tr2} T. Tradler, ``The BV algebra on Hochschild Cohomology
induced by Infinity-inner products'' \texttt{math.QA/0210150}

\bibitem[TZ]{TZ} T. Tradler, M. Zeinalian, ``Poincare Duality at the
Chain Level, and a BV Structure on the Homology of the Free Loops
Space of a Simply Connected Poincare Duality Space''
\texttt{math.AT/0309455}

\bibitem[W]{W} S. Wilson, SUNY Stony Brook Ph.D. Thesis

\end{thebibliography}
\end{document}